\documentclass[12pt,leqno]{article}
\usepackage{amsmath, amssymb}

\newcommand{\sect}[1]{\setcounter{equation}{0}\section{#1}}

\def\epsilon{\varepsilon}

\begin{document}

\LARGE \noindent 
{\bf Finite-dimensional complex manifolds on commutative Banach algebras 
and continuous families of compact complex manifolds}

\large

\vspace*{0.4em}

\hfill Hiroki Yagisita (Kyoto Sangyo University)

\vspace*{1.2em}

\normalsize

Abstract: 

An $n$-dimensional complex manifold is a manifold by biholomorphic mappings 
between open sets of the finite direct product ${\mathbb C}^n$ of the complex number field. 
On the other hand, when $A$ is a commutative Banach algebra, 
Lorch gave a definition that an $A$-valued function 
on an open set of $A$ is holomorphic.  
The definition of a holomorphic function by Lorch 
can be straightforwardly generalized 
to an $A$-valued function on an open set of the finite direct product $A^n$. 
Therefore, a manifold modeled on $A^n$ 
(an $n$-dimensional $A$-manifold) is easily defined. 
However, in my opinion, it seems that so many nontrivial examples 
were not known (including the case of $n=1$, that is, Riemann surfaces). 

By the way, if $X$ is a compact Hausdorff space, 
then the algebra $C(X)$ of all complex valued continuous functions on $X$ 
is the most basic example of a commutative Banach algebra 
(furthermore, a commutative $C^*$-algebra). 
In this note, we see that if the set 
of all continuous cross sections 
of a continuous family $M$ of compact complex manifolds 
(a topological deformation $M$ of compact complex analytic structures) on $X$ 
is denoted by $\Gamma(M)$, 
then the structure of a $C(X)$-manifold modeled on the $C(X)$-modules  
of all continuous cross sections of complex vector bundles on $X$ 
is introduced into $\Gamma(M)$. 
Therefore, especially, if $X$ is contractible, 
then $\Gamma(M)$ is a finite-dimensional $C(X)$-manifold.

\vfill

Japanese: 

https://www.researchgate.net/publication/327236804

\newpage

$\, $

\vfill 

%\noindent 
%Mathematics Subject Classification. Primary 32Q99; Secondary 46J99, 58B99.


\newpage

\normalsize

\sect{Compact continuous families of complex manifolds} 
\ \ \ \ \ In this and the next sections, we introduce two basic concepts of this note.

\vspace*{0.4em}

\noindent 
{\bf Definition 1.1} (Canonical real coordinate system) : 

The mapping that maps a complex vector $z=x+\sqrt{-1}y \, \in \, {\mathbb C}^n$ 
to its real and imaginary pair $(x,y) \, \in \, {\mathbb R}^{2n}$, i.e.,  
$$( \, x_1+\sqrt{-1}y_1, \, x_2+\sqrt{-1}y_2, \, \cdots, \, x_n+\sqrt{-1}y_n \, )$$
$$ \mapsto \ ((x_1,x_2,\cdots,x_n), \, (y_1,y_2,\cdots,y_n))$$
is called the canonical real coordinate system. 
\hfill ---

\vspace*{0.2em}

Hereafter, when we consider an open set of ${\mathbb C}^n$ 
to be a $C^\infty$-manifold, we use the canonical real coordinate system 
as its local coordinate system. 

\vspace*{0.4em}

\noindent 
{\bf Definition 1.2} (Open polydisk) : 

For $r>0$, let 
$$D_r^n \ := \ \{ \ (z_1,z_2,\cdots,z_n) \in {\mathbb C}^n \ | \ \max_k|z_k| \, < \, r \ \}$$
$$= \ \{ \ ((x_1,x_2,\cdots,x_n), \, (y_1,y_2,\cdots,y_n)) \ \in {\mathbb R}^{2n} \ |$$
$$\max_k\sqrt{{x_k}^2+{y_k}^2} \, < \, r \ \}.$$
Let $D^n:=D_1^n$. $D^n$ is called the unit open polydisk.
\hfill ---

\noindent 
{\bf Definition 1.3} (Compact continuous family of complex manifolds) : 

$M \, := \, (M,X,\pi,S)$ is said to be a compact continuous family 
(of $n$-dimensional complex manifolds), if it satisfies the followings.

(1) \ \ \ $M$ and $X$ are compact Hausdorff spaces. 
$\pi$ is a continuous surjection form $M$ to $X$.

(2) \ \ \ $S \, (\not=\emptyset)$ is a set. 
For any $\varphi\in S$, $\varphi$ is a map 
from an open set $M_\varphi \, (\not=\emptyset)$ of $M$ 
to the unit open polydisk $D^n$. 
For $\varphi\in S$, let 
$$\pi |_\varphi \ := \ \pi_{\upharpoonright{M_\varphi}}.$$
For any $\varphi\in S$, there 
exists an open set $X_\varphi \, (\not=\emptyset)$ of $X$ 
such that $(\varphi, \pi |_\varphi)$ is a homeomorphism 
from $M_\varphi$ to $D^n \times X_\varphi$. 

(3) \ \ \ For $\varphi\in S$ and $t\in X_\varphi$, let 
$$M_{\varphi,t} \ := \ M_\varphi\cap\pi^{-1}(\{t\}),$$
$$\varphi |_t \ := \ \varphi_{\upharpoonright M_{\varphi,t}}.$$

\newpage \noindent 
For any $\varphi_1, \varphi_2 \in S$ and  
$t\in X_{\varphi_1}\cap X_{\varphi_2}$, 
the coordinate transformation 
$${\varphi_2} |_t \, \circ \, {\varphi_1} |_t^{-1} \ : \ 
\varphi_1(M_{\varphi_1,t}\cap M_{\varphi_2,t})
\, \rightarrow \, 
\varphi_2(M_{\varphi_1,t}\cap M_{\varphi_2,t})$$
is holomorphic. 

(4) \ \ \ $M \, = \, \cup_{\varphi\in S} \, M_\varphi$ holds. 

\noindent 
{\sf Remark} : 

(1) \ \ \ Each fiber of a compact continuous family 
is a compact complex manifold. 

(2) \ \ \ Suppose that $(M,B,\pi)$ is a complex analytic family  
of compact complex manifolds. Then, 
for any compact subset $X \, (\not=\emptyset)$ of $B$, 
$\pi^{-1}(X)$ is a compact continuous family. 
\hfill ---

\vspace*{0.2em}

Unlike complex analytic families, 
it seems that continuous families 
are not subjects of strong interest. 
Therefore, we have to confirm very primitive facts 
of compact continuous families on our own. 

\vspace*{0.2em}

\noindent 
{\bf Definition 1.4} (Local trivialization coordinate neighborhood) : 

When $M \, := \, (M,X,\pi,S)$ is a compact continuous family, 
$M$ is called the total space, $X$ is called the base space, 
$\pi$ is called the projection 
and $S$ is called the system 
of local trivialization coordinate neighborhoods. 
$\mbox{\ \ \ }$ \hfill ---

\vspace*{0.2em}

Since the total space of a compact continuous family is compact, 
the system of local trivialization coordinate neighborhoods 
can be assumed to be a finite set. 

\vspace*{0.2em}

\noindent 
{\bf Definition 1.5} (Finite continuous family) : 

A compact continuous family whose system 
of local trivialization coordinate neighborhoods 
is a finite set is called a finite continuous family. 
\hfill ---

\noindent 
{\bf Definition 1.6} (Compact coordinate neighborhood) : 

Let $(M,X,\pi,S)$ be a finite continuous family. Then, 
$(\, \{\, (M^\prime_\varphi,X^\prime_\varphi)\, \}_{\varphi\in S}, \, r_0 \, )$ 
is said to be a system of compact coordinate neighborhoods of $M$, 
if it satisfies the followings.

(1) \ \ \ $M^\prime_\varphi$ is an open set of $M$, 
$X^\prime_\varphi$ is an open set of $X$ and 
$$\overline{M^\prime_\varphi}\subset M_\varphi, \ \ \ 
\overline{X^\prime_\varphi}\subset X_\varphi$$
hold. 

(2) \ \ \ $r_0\in (0,1)$ holds.

(3) \ \ \ The homeomorphism $(\varphi, \pi |_\varphi)$ 
from $M_\varphi$ to $D^n \times X_\varphi$ 
maps $M^\prime_\varphi$ onto $D_{r_0}^n \times X^\prime_\varphi$. 
That is, 
$$M^\prime_\varphi \
= \ (\varphi, \pi |_\varphi)^{-1}
\, (D_{r_0}^n \times X^\prime_\varphi)$$
holds. 

\newpage 

(4) \ \ \ $M \, = \, \cup_{\varphi\in S} \, M^\prime_\varphi$ holds. 

\noindent 
{\sf Remark} :
In the above definition, 
neither $M^\prime_\varphi\not=\emptyset$ 
nor $X^\prime_\varphi\not=\emptyset$ 
is required. 
\hfill ---

\noindent 
{\bf Definition 1.7} (Bump system) : 

Let $(M,X,\pi,S)$ be a finite continuous family. Then, 
$(\, \{\, (\rho_\varphi,M^\prime_\varphi,X^\prime_\varphi)\, \}_{\varphi\in S}, \, r_0 \, )$ 
is said to be a bump system of $M$, 
if it satisfies the followings. 

(1) \ \ \ $(\, \{\, (M^\prime_\varphi,X^\prime_\varphi)\, \}_{\varphi\in S}, \, r_0 \, )$ 
is a system of compact coordinate neighborhoods of $M$. 

(2) \ \ \ $\rho_\varphi$ is a continuous function 
from $M$ to $[0,+\infty)$ and 
$${\rm supp}(\rho_\varphi) \ \subset \ M^\prime_\varphi$$
holds. 

(3) \ \ \ The function 
$$\rho_\varphi \, \circ \, (\varphi,\pi |_\varphi)^{-1}
\ : \ \ \ D^n \times X_\varphi \ \rightarrow \ [0,+\infty)$$
is $C^\infty$-class with respect to the variable 
$(x,y) \, \in \, D^n \, (\subset {\mathbb R}^{2n})$. 
For any multi-indexes $\alpha$ and $\beta$, the function 
$$ {\partial_x}^\alpha \, {\partial_y}^\beta 
\ ( \, \rho_\varphi \, \circ \, (\varphi,\pi |_\varphi)^{-1} \, )
\ : \ \ \ D^n \times X_\varphi \ \rightarrow \ {\mathbb R}
$$ is continuous. 

(4) \ \ \ For any $p\in M$, there exists $\varphi \in S$ 
such that $$\rho_\varphi(p)\not=0$$ holds. 

\noindent 
{\sf Remark} : In the above definition, 
$\sum_{\varphi\in S} \, \rho_\varphi \, = \, 1$
is not required. 
\hfill ---

\noindent 
{\bf Lemma 1.8} (Existence of a bump system) : 

Suppose that $M$ is a finite continuous family. 
Then, there exists a bump system of $M$. 

\noindent 
{\sf Proof} : We give it in Section 5. 
\hfill 
$\blacksquare$ 

\noindent 
{\bf Definition 1.9} (Continuous family with bumps) : 

$M:=(M,X,\pi,S,\{\, (\rho_\varphi,M^\prime_\varphi,X^\prime_\varphi)\, \}_{\varphi\in S},r_0)$ 
is said to be a continuous family with bumps, 
if $(M,X,\pi,S)$ is a finite continuous family 
and 
$(\{\, (\rho_\varphi,M^\prime_\varphi,X^\prime_\varphi)\, \}_{\varphi\in S}, r_0)$ 
is a bump system of $M$. 
\hfill ---

\vspace*{0.2em}

Hereafter, the continuous family $M \, 
(=(M,X,\pi,S,\{\, (\rho_\varphi,M^\prime_\varphi,X^\prime_\varphi)\, \}_{\varphi\in S},r_0))$ 
with bumps is fixed to be one. 

\vspace*{0.4em}

\newpage \noindent 
{\bf Definition 1.10} (Real tangent bundle) : 

Let $N$ be an $n$-dimensional complex manifold. 
Then, for $q \in N$, 
its real tangent space $(T_{\mathbb R})_q(N)$ is a complex linear space 
by the almost complex structure of $N$. That is, 
when $z=x+\sqrt{-1}y$ is a holomorphic local coordinate system of $N$, set 
$$\sqrt{-1} \, (\frac{\partial}{\partial x_k})_q 
\ := \ (\frac{\partial}{\partial y_k})_q, \ \ \ 
\sqrt{-1} \, (\frac{\partial}{\partial y_k})_q 
\ := \ - \, (\frac{\partial}{\partial x_k})_q.$$
In addition, set 
$$T_{\mathbb R}(N) \ := \ \cup_{q\in N} \ (T_{\mathbb R})_q(N).$$ 
$T_{\mathbb R}(N)$ is a holomorphic vector bundle on $N$. 

\noindent 
{\sf Remark} : The real tangent bundle $T_{\mathbb R}(N)$ 
is holomorphically isomorphic to 
the holomorphic tangent bundle $T^\prime(N)$ by 
$$(\frac{\partial}{\partial x_k})_q \, \in \, T_{\mathbb R}(N)
\ \ \ \mapsto \ \ \ 
(\frac{\partial}{\partial z_k})_q
:=\frac{1}{2}((\frac{\partial}{\partial x_k})_q
-\sqrt{-1}(\frac{\partial}{\partial y_k})_q)
\, \in \, T^\prime(N).$$
\hfill ---

\noindent 
{\bf Lemma 1.11} : 

Let $T$ be a topological space. Let $U$ be an open set of ${\mathbb C}^n\times T$. 
Suppose that $f$ is a complex valued continuous function on $U$ 
and for any $t\in T$, the map $z \ \mapsto \ f(z,t)$ is holomorphic. Then, 
for any multi-index $\gamma$, 
$${\partial_z}^\gamma \, f \ : \ U \, \rightarrow \, {\mathbb C}$$
is continuous. 

\noindent 
{\sf Proof} : We give it in Section 6. 
\hfill 
$\blacksquare$

\noindent 
{\bf Definition 1.12} (Trivial Hermitian metric) :  

Let $\varphi \in S$. For $p\in M_\varphi$, define 
a complex inner product on 

\noindent 
$(T_{\mathbb R})_p \, (\pi^{-1}(\{\pi(p)\}))$ as 
$$\langle \, \dot{p}_1, \, \dot{p}_2 \, \rangle_{\varphi,p}
\ := \ 
\langle \, (D(\varphi |_{\pi(p)}))_p \, (\dot{p}_1), 
\, (D(\varphi |_{\pi(p)}))_p \, (\dot{p}_2) \, 
\rangle_{{\mathbb C}^n}.$$
Here, $D(\varphi |_{\pi(p)})$ 
is the differential of the holomorphic local coordinate system 
$$\varphi |_{\pi(p)}
\ := \ \varphi_{\upharpoonright M_{\varphi, \pi(p)}} 
\ = \ \varphi_{\upharpoonright M_\varphi\cap\pi^{-1}(\{\pi(p)\})}$$
of $\pi^{-1}(\{\pi(p)\})$
and $\langle \cdot, \cdot \rangle_{{\mathbb C}^n}$ 
is the canonical inner product (the canonical Hermitian metric) 
on ${\mathbb C}^n$. 
Also, $\|\cdot\|_{\varphi,p}$ denotes the norm 
by this complex inner product. That is, let 
$$\|\dot{p}\|_{\varphi,p} \ := \ 
\sqrt{\langle \dot{p}, \dot{p} \rangle_{\varphi, p}}.$$
For $t\in X_\varphi$, by the complex inner products, 
the complex manifold 
$$M_{\varphi,t} \ := \ M_\varphi\cap\pi^{-1}(\{t\})$$
is a Hermitian manifold that is isomorphic to $D^n$. 
In addition, the topological space 
$$T_{\mathbb R} \, (M_\varphi)
\ := \ \cup_{p\in M_\varphi} \, (T_{\mathbb R})_p \, (\pi^{-1}(\{\pi(p)\}))$$
is a continuous Hermitian vector bundle on the topological space $M_\varphi$. 
\hfill ---

\noindent 
{\bf Definition 1.13} (Hermitian real tangent bundle) :  

For $p\in M$, define a complex inner product on 
$(T_{\mathbb R})_p \, (\pi^{-1}(\{\pi(p)\}))$ as 
$$\langle \, \dot{p}_1, \, \dot{p}_2 \, \rangle_p
\ := \ \sum_{\varphi \in S} 
\, \rho_\varphi(p)
\, \langle \, \dot{p}_1, \, \dot{p}_2 \, \rangle_{\varphi,p}.$$
Also, $\|\cdot\|_p$ denotes the norm by this complex inner product. 
That is, let 
$$\|\dot{p}\|_p \ := \ 
\sqrt{\langle \dot{p}, \dot{p} \rangle_p}.$$ 
For $t\in X$, by the complex inner products, 
the fiber $\pi^{-1}(\{t\})$ of the continuous family $M$ 
is a Hermitian manifold. In addition, the topological space 
$$T_{\mathbb R} \, (M)
\ := \ \cup_{p\in M} \, (T_{\mathbb R})_p \, (\pi^{-1}(\{\pi(p)\}))$$
is a continuous Hermitian vector bundle 
on the total space $M$ of the continuous family. 
The projection $T_{\mathbb R}(M) \, \rightarrow \, M$ 
of the vector bundle is denoted by $\varpi$. That is, 
for $\dot{p} \, \in \, T_{\mathbb R}(M)$, 
$$\varpi(\dot{p}) \, \in \, M, \ \ \ \ \ \ 
\dot{p} \, \in \, (T_{\mathbb R})_{\varpi(\dot{p})}
(\pi^{-1}(\{\pi(\varpi(\dot{p}))\}))$$
hold. $\varpi(\dot{p})$ is called the base point 
of a real tangent vector $\dot{p}$.

\noindent 
{\sf Remark} : In the definition of $T_{\mathbb R}(M)$, Lemma 1.11 was used. 
That is, according to Lemma 1.11, for any local trivialization coordinate systems 
$\psi_1$ and $\psi_2$ of $M$ that are compatible with each other, 
the local trivialization coordinate system of $T_{\mathbb R}(M)$ determined by $\psi_1$ 
and the one determined by $\psi_2$ are compatible with each other. 
\hfill ---

\noindent 
{\bf Lemma 1.14} : 

Let $\varphi\in S$. Let $K$ be a compact subset of $M_\varphi$. 
Then, there exist an open set $U$ of $M_\varphi$ and $C>0$ 
such that the followings hold. 

\newpage 

(1) \ \ \ $K\subset U$ holds. 

(2) \ \ \ 
For any $p\in U$ and $\dot{p}\in (T_{\mathbb R})_p(\pi^{-1}(\{\pi(p)\}))$, 
$$\|\dot{p}\|_{\varphi,p}
\ \leq \ C \, \|\dot{p}\|_p,
\ \ \ \ \ \ \|\dot{p}\|_p
\ \leq \ C \, \|\dot{p}\|_{\varphi,p}$$
hold. 

\noindent 
{\sf Proof} : 
Both the Hermitian metrics are continuous. So, it follows. 
\hfill 
$\blacksquare$

\noindent 
{\bf Definition 1.15} (Distance on a fiber) : 

Let $t\in X$. For a piecewise $C^\infty$-curve 
$c\, : \, [a,b]\rightarrow \pi^{-1}(\{t\})$, let 
$$L_t(c) \ := \ \int_a^b 
\, \|\frac{d}{ds}(c(s))\|_{c(s)} 
\, ds.$$
For $p, q \in \pi^{-1}(\{t\})$, let 
$$d_t(p,q)$$
$$:= \ \inf \, ( \, \{ \, L_t(c) \, | \, 
\mbox{$c$ is a piecewise $C^\infty$-curve}$$
$$\mbox{with the start point $p$ 
and the end point $q$} 
\, \} \cup \{1\} \, ).
$$
By $d_t$, $\pi^{-1}(\{t\})$ is a metric space.

\noindent 
{\sf Remark} : 
${\rm Re} \, (\langle z,w\rangle_{{\mathbb C}^n}) \, = \, 
\langle z,w\rangle_{{\mathbb R}^{2n}}$ and especially,  
$\|z\|_{{\mathbb C}^n} \, = \, \|z\|_{{\mathbb R}^{2n}}$ hold. 
Therefore, the norm defined by a complex inner product 
is the norm defined by a real one. 
\hfill ---

\noindent 
{\bf Lemma 1.16} : 

Let $\varphi\in S$. Let $N$ be an open set of $M_\varphi$. 
Let $K$ be a compact subset of $N$. 
Then, there exist an open set $U$ of $N$, $\delta>0$ and $C>0$ 
such that the followings hold. 

(1) \ \ \ $K\subset U$ holds. 

(2) \ \ \ 
Suppose $p \, \in \, U$ and $q \, \in \, \pi^{-1}(\{\pi(p)\})$. Then, 
$$d_{\pi(p)}(q,p) \, < \, \delta$$
$$\Longrightarrow$$ 
$$q \, \in \, N, 
\ \ \ \|\varphi(q)-\varphi(p)\|_{{\mathbb C}^n} \, \leq \, C \, d_{\pi(p)}(q,p)$$
holds.

\noindent 
{\sf Proof} : We give it in Section 7. 
\hfill 
$\blacksquare$

\noindent 
{\bf Remark on Lemma 1.16} :  

Let $\varphi\in S$. Let $N$ be an open set of $M_\varphi$. 
Let $K$ be a compact subset of $N$. Then, 
there exist an open set $U$ of $N$, $\delta>0$ and $C>0$ 
such that the followings hold. 

(1) \ \ \ $K\subset U$ holds. 

(2) \ \ \ 
Suppose that $p \, \in \, U$ and $q \, \in \, M_{\varphi,\pi(p)} 
\, (:=M_\varphi\cap\pi^{-1}(\{\pi(p)\}))$. Then, 
$$\|\varphi(q)-\varphi(p)\|_{{\mathbb C}^n} \, < \, \delta$$
$$\Longrightarrow$$ 
$$q \, \in \, N, 
\ \ \ d_{\pi(p)}(q,p) \, \leq \, C \, \|\varphi(q)-\varphi(p)\|_{{\mathbb C}^n}$$
holds. 

\noindent 
{\sf Proof} : 
We give it in Section 8. 
(However, this remark is not used.) 
\hfill 
$\blacksquare$

\noindent 
{\bf Definition 1.17} (Continuous section) :  

Let $T\, (\not=\emptyset)$ be a subset of $X$. 
$u$ is said to be a continuous section of $M$ on $T$, 
if $u$ is a continuous map from $T$ to $M$ 
and for any $t\in T$, $\pi(u(t))\, = \, t$ holds. 
The set of all continuous sections of $M$ on $T$ 
is denoted by $\Gamma(M|T)$. 
\mbox{ \, \, \, } \hfill ---

\noindent 
{\bf Definition 1.18} (Distance between continuous sections) : 

Let $T\, (\not=\emptyset)$ be a subset of $X$. 
For $u, v \in \Gamma(M|T)$, let 
$$d_T(u,v) \ := \ \sup_{t\in T} \, d_t(u(t),v(t)).$$
By $d_T$, $\Gamma(M|T)$ is a metric space. 
\hfill ---

\newpage

\sect{Manifolds on commutative Banach algebras}
\noindent 
{\bf Definition 2.1} (Commutative algebra) : 

$A$ is said to be a commutative algebra, 
if $A$ is a commutative ring with the multiplicative unit $1_A \, (\not=0_A)$, 
it is a complex linear space and it satisfies 
$$(c1_A)f=cf \ \ \ \ \ \ ( \, c\in \mathbb C, \, \, \, f\in A \, ).$$ 

\noindent 
{\sf Remark} : When $A$ is a commutative algebra, ${\mathbb C} \subset A$ holds. 
\hfill ---

\noindent 
{\bf Definition 2.2} (Commutative Banach algebra) : 

$A$ is said to be a commutative Banach algebra, 
if $A$ is a commutative algebra, it is a complex Banach space 
and it satisfies 
$$\|fg\| \, \leq \, \|f\| \, \|g\| \ \ \ \ \ \ ( \, f, g \, \in \, A \, ),$$
$$\|1_A\|=1.$$
\hfill ---

\noindent 
{\bf Example 2.3} : 

Let $T \, (\not=\emptyset)$ be a topological space. 
Let $C_b(T)$ be the set of all complex valued bounded continuous functions on $T$. 
Let $$\|f\| \, := \, \sup_{t\in T} \, |f(t)| \ \ \ \ \ \ ( \, f\in C_b(T) \, ).$$
$C_b(T)$ is a commutative Banach algebra. 
\hfill ---

\noindent 
{\bf Definition 2.4} (Banach module) : 

Let $A$ be a commutative Banach algebra. 
$X$ is said to be a Banach $A$-module, 
if $X$ is an $A$-module, it is a complex Banach space 
and it satisfies 
$$(c1_A)u \, = \, cu \ \ \ \ \ \ ( \, c\in {\mathbb C}, \, \, \, u\in X \, ),$$
$$\|fu\|_X \, \leq \, \|f\|_A \, \|u\|_X \ \ \ \ \ \ ( \, f\in A, \, \, \, u\in X \, ).$$
\hfill ---

\noindent 
{\bf Example 2.5} : 

Let $A$ be a commutative Banach algebra. 
Let $$\|(f_1,f_2,\cdots,f_n)\| \ := \
\max_k \, \|f_k\|_A \ \ \ \ \ \ ( \, f_1,f_2,\cdots,f_n \, \in \, A \, ).$$
The finite direct product $A^n$ is a Banach $A$-module. 
\hfill ---

\newpage 

\noindent 
{\bf Example 2.6} : 

Let $T \, (\not=\emptyset)$ be a topological space. 
Let $X$ be a complex Banach space. 
Let $C(T;X)$ be the set of all $X$-valued continuous functions on $T$. 
Let 
$$\|u\| \ := \ \sup_{t\in T} \, \|u(t)\|_X 
\ \ \ \ \ \ ( \, u\in C(T;X) \, ),$$
$$C_b(T;X) \ := \ \{ \ u\in C(T;X) \ | \ \|u\|<+\infty \ \}.$$
$C_b(T;X)$ is a Banach $C_b(T)$-module. 
\hfill ---

\noindent 
{\bf Example 2.7} : 

Let $T \, (\not=\emptyset)$ be a topological space. 
Let $E$ be a continuous Hermitian vector bundle on $T$. 
Let $\Gamma(E)$ be the set of all continuous sections of $E$ on $T$. 
Let 
$$\|u\| \ := \ \sup_{t\in T} \, \|u(t)\|_t 
\ = \ \sup_{t\in T} \, \sqrt{ \, \langle \, u(t), \, u(t) \, \rangle_t \, }
\ \ \ \ \ \ ( \, u\in \Gamma(E) \, ),$$
$$\Gamma_b(E) \ := \ \{ \ u\in\Gamma(E) \ | \ \|u\|<+\infty \ \}.$$ 
$\Gamma_b(E)$ is a Banach $C_b(T)$-module. 
\hfill ---

\noindent 
{\bf Definition 2.8} (Linear mapping) : 

Let $A$ be a ring. Let $X$ and $Y$ be $A$-modules. 
A mapping $F:X\rightarrow Y$ is said to be $A$-linear, 
if it satisfies 
$$F(u+v) \, = \, F(u)+F(v) \ \ \ \ \ \ (\, u, v \, \in \, X \, ),$$
$$F(fu) \, = \, fF(u) \ \ \ \ \ \ ( \, f\in A, \, \, \, u\in X \, ).$$

\noindent 
{\sf Remark} : 
Let $X$ and $Y$ be Banach $A$-modules. 
If a mapping $F:X\rightarrow Y$ is $A$-linear, 
then $F$ is complex linear. 
\hfill ---

\noindent 
{\bf Example 2.9} : 

Let $T \, (\not=\emptyset)$ be a topological space. 
Let $A$ be a continuous complex linear map 
from a complex Banach space $X$ to a complex Banach space $Y$. 
If $\tilde{A}$ is defined as 
$$(\tilde{A}(u))(t) \ := \ A(u(t)) \ \ \ \ \ \ (u\in C_b(T;X), \ t\in T),$$
then $\tilde{A}$ is a continuous $C_b(T)$-linear map 
from $C_b(T;X)$ to $C_b(T;Y)$. \hfill ---

\noindent 
{\bf Definition 2.10} (Frechet derivative) :  

Let $X$ and $Y$ be real Banach spaces. 
Let $f$ be a map from an open set $U$ of $X$ to $Y$.

(1) \ \ \ Let $p\in U$. 
Let $A\, :X\rightarrow Y$ be a continuous real linear map. 
$A$ is said to be the Frechet derivative of $f$ at the point $p$, 
if for any $\varepsilon>0$, there exists $\delta>0$ 
such that 
$$\|h\|_X \, \leq \, \delta \ \ \ \Longrightarrow \ \ \ 
\|f(p+h)-(f(p)+Ah)\|_Y \, \leq \, \varepsilon \|h\|_X$$
holds. 

(2) \ \ \ $f$ is said to be Frechet differentiable (on $U$), 
if for any $p\in U$, there exists the Frechet derivative of $f$ at the point $p$. 

\noindent 
{\sf Remark} : 
If $A$ and $B$ are the Frechet derivative of $f$ at a point $p$, 
then $A=B$ holds. The Frechet derivative of $f$ at a point $p$ 
is denoted by $f^\prime(p)$, $(Df)(p)$, $(Df)_p$ and so on. 
\hfill ---

\noindent 
{\bf Definition 2.11} (Complex differentiable) :  

Let $X$ and $Y$ be complex Banach spaces. 
Let $f$ be a map from an open set $U$ of $X$ to $Y$. 
$f$ is said to be complex differentiable (on $U$), 
if $f$ is Frechet differentiable 
and for any $p\in U$, the Frechet derivative $(Df)_p$ 
is complex linear.

\noindent 
{\sf Remark} (Complex power expansion) :  

Suppose that $f$ is complex differentiable. 
Then, $f$ is complex analytic (that is, 
there exists the complex power expansion 
in a neighborhood of each point).  
In particular, $Df$ is complex differentiable. 
([10] Theorem 14.7)
\hfill ---

\vspace*{0.4em}

Let us generalize the definition of an $A$-differentiable function  
from an open set of a commutative Banach algebra $A$ 
to $A$ introduced by Lorch [5] 
to a mapping from an open set of a Banach $A$-module 
to a Banach $A$-module. 

\vspace*{0.2em}

\noindent 
{\bf Definition 2.12} (Differentiable mapping) : 

Let $X$ and $Y$ be Banach $A$-modules. 
Let $f$ be a mapping from an open set $U$ of $X$ to $Y$. 
$f$ is said to be $A$-differentiable (on $U$), 
if $f$ is Frechet differentiable 
and for any $p\in U$, the Frechet derivative $(Df)_p$ 
is $A$-linear. 

\noindent 
{\sf Remark} : 
If one is $A$-differentiable, then it is complex differentiable. 
\hfill ---

\vspace*{0.2em}

The following theorem is known 
about the Lorch differential of an $A$-valued function 
on an open set of $A$ ([3] \S 3.19 and \S 26.4). 
However, we do not use this theorem. 
However, with this, the identity of $C_b(T)$-differentiable mappings 
from open sets of $C_b(T)$ to $C_b(T)$ is almost obvious. 
Rather, it would be more accurate 
to say that this theorem is a generalization of that fact. 

\vspace*{0.2em}

\noindent 
{\bf Theorem 2.13} (Taylor expansion) :  

Let $A$ be a commutative Banach algebra. 
Let $c\in A$ and $r\in (0,+\infty]$. 
Let $D_r(c) \, := \, \{ \, x\in A \, | \, \|x-c\|<r \, \}$. 

(1) \ \ \ Let $\{a_n\}_{n=0}^\infty$ be a sequence of $A$. 
Suppose $r\, =\, \liminf_{n\rightarrow \infty}\, (\, \|a_n\|^{-\frac{1}{n}}\, )$. 
Then, 
$$x \ \ \ \mapsto \ \ \ \sum_{n=0}^\infty \ a_n \, (x-c)^n$$ 
is an $A$-differentiable mapping from $D_r(c)$ to $A$. 

(2) \ \ \ Suppose that $f$ is an $A$-differentiable mapping from $D_r(c)$ to $A$. 
Then, there uniquely exists a sequence $\{a_n\}_{n=0}^\infty$ of $A$ 
such that 
$$f(x) \ = \ \sum_{n=0}^\infty \ a_n \, (x-c)^n$$ 
holds. Further, 
$r \, \leq \, \liminf_{n\rightarrow \infty} \, (\|a_n\|^{-\frac{1}{n}})$ 
holds. 
\hfill ---

\noindent 
{\bf Definition 2.14} (Manifold on a commutative Banach algebra) : 

Let $M$ be a Hausdorff space. Let $S$ be a set. 
Let $A$ be a commutative Banach algebra. 
$M$ is said to be a complex manifold on the commutative Banach algebra $A$ 
(or an $A$-manifold) with the system $S$ of coordinate neighborhoods, 
if the followings hold. 

For any $\varphi \in S$, there exists a Banach $A$-module $X$ 
such that $\varphi$ is a homeomorphism from an open set of $M$
to an open set of $X$. 
For any $p\in M$, there exists $\varphi \in S$ 
such that $p$ belongs to the domain of $\varphi$. 
For any $\varphi_1, \varphi_2\in S$, the coordinate transformation 
$$\varphi_2 \, \circ \, \varphi_1^{-1} \ : \ \varphi_1(U_1\cap U_2) \ 
\rightarrow \ \varphi_2(U_1\cap U_2)$$ 
is $A$-differentiable. Here, $U_1$ and $U_2$ 
are the domains of $\varphi_1$ and $\varphi_2$, respectively. 
$\mbox{}$ \hfill ---

\noindent 
{\sf Remark} (Riemann sphere) : 
In [1], the Riemann spheres on commutative Banach algebras are considered. 
\hfill ---

\noindent 
{\sf Remark} (Holomorphic vector bundle on a commutative Banach algebra) :  

(1) \ \ \ 
Let $M_1$ be an $A$-manifold with a system $S_1$ of coordinate neighborhoods. 
Let $M_2$ be an $A$-manifold with a system $S_2$ of coordinate neighborhoods. 
Let $f \, : \, M_1\rightarrow M_2$ be a continuous mapping. 
$f$ is said to be $A$-holomorphic, 
if for any $\varphi_1 \in S_1$ and $\varphi_2 \in S_2$, the mapping 
$$\varphi_2 \, \circ \, f \, \circ \, \varphi_1^{-1} 
\ : \ \varphi_1(U_1\cap f^{-1}(U_2)) \ 
\rightarrow \ X_2$$
is $A$-differentiable. 
Here, $\varphi_1$ is a mapping from $U_1$ to a Banach $A$-module $X_1$ 
and $\varphi_2$ is a mapping from $U_2$ to a Banach $A$-module $X_2$.

(2) \ \ \ 
$(E,M,\pi,\{(U_\lambda, X_\lambda, \varphi_\lambda)\}_{\lambda\in\Lambda})$ 
is said to be an $A$-holomorphic vector bundle, 
if it satisfies the followings. 

$E$ and $M$ are $A$-manifolds. 
$\pi$ is an $A$-holomorphic surjection from $E$ to $M$. 
Each $U_\lambda$ is an open set of $M$. 
$M \, = \, \cup_{\lambda\in\Lambda} U_\lambda$ holds. 
Each $\varphi_\lambda$ is a map 
from $\pi^{-1}(U_\lambda)$ 
to a Banach $A$-module $X_\lambda$. For $\lambda\in\Lambda$, let 
$$\pi |_\lambda \ := \ \pi_{\upharpoonright \pi^{-1}(U_\lambda)}.$$
Each map 
$$(\varphi_\lambda,\pi |_\lambda) \ : \ 
\pi^{-1}(U_\lambda) \, \rightarrow \, X_\lambda \times U_\lambda$$
is an $A$-biholomorphic map. For $\lambda\in\Lambda$ and $p\in U_\lambda$, let 
$${\varphi_\lambda} |_p 
\ := \ {\varphi_\lambda}_{\upharpoonright \pi^{-1}(\{p\})}.$$ 
For any $\lambda_1, \lambda_2 \, \in \, \Lambda$ 
and $p \, \in \, U_{\lambda_1}\cap U_{\lambda_2}$, 
the coordinate transformation 
$${\varphi_{\lambda_2}} |_p \, \circ 
\, {\varphi_{\lambda_1}} |_p^{-1} \ : \ 
X_{\lambda_1} \, \rightarrow \, X_{\lambda_2} $$
is $A$-linear. 
\hfill ---

\newpage

\sect{The main result} 
\ \ \ \ \ First of all, recall that the continuous family 
$M \, (=(M,X,\pi,S,\{\, (\rho_\varphi,M^\prime_\varphi,X^\prime_\varphi)\, \}_{\varphi\in S},r_0))$ 
with bumps was fixed to be one. 

\vspace*{0.4em}

\noindent 
{\bf Definition 3.1} (Hermitian vector bundle induced by a continuous section) : 

Let $T\, (\not=\emptyset)$ be a subset of $X$. 
Let $u\in \Gamma(M|T)$. Let 
$$u^* \, (T_{\mathbb R}(M))
\ := \ T_{\mathbb R}(M) | u(T)$$
$$= \ \cup_{\, p \in u(T)} 
\, (T_{\mathbb R})_p(\pi^{-1}(\{\pi(p)\})).$$ 
Because $u$ is a homeomorphism from $T$ to $u(T)$, 
$u^* \, (T_{\mathbb R}(M))$ is a continuous Hermitian vector bundle 
on the topological space $T$. 
\hfill ---

\noindent 
{\bf Definition 3.2} (Norm of a bounded continuous section) :  

Let $T \, (\not=\emptyset)$ be a subset of $X$. 
Let $u\in \Gamma(M|T)$. Then, 
the norm of a bounded continuous section 
of $u^* \, (T_{\mathbb R}(M))$ on $T$ 
is denoted by $\|\cdot\|_u$. That is, for $\dot{u} 
\, \in \, \Gamma_b \, ( \, u^* \, (T_{\mathbb R}(M)) \, )$, let 
$$\|\dot{u}\|_u \ := \ \sup_{t\in T} \, \|\dot{u}(t)\|_{u(t)}.$$
$\Gamma_b \, ( \, u^* \, (T_{\mathbb R}(M)) \, )$ 
is a Banach $C_b(T)$-module. 
\hfill ---

\noindent 
{\bf Definition 3.3} (Holomorphic normal coordinate neighborhood) : 

Let $T\, (\not=\emptyset)$ be a subset of $X$. 
Let $u \, \in \, \Gamma(M|T)$. Then, 
$(U, \Psi)$ is said to be a holomorphic normal coordinate neighborhood of $u$, 
if it satisfies the followings. 

(1) \ \ \ $U$ is an open set of $\pi^{-1}(T)$. 
For any $t\in T$, 
$$u(t) \ \in \ U_t \ ( \, := \, U\cap \pi^{-1}(\{t\}) \, )$$
holds. 

(2) \ \ \ For $p\in M$, let 
$$B_p(T_{\mathbb R}) 
\ : \, = \ \{\ \dot{p} \, \in \, (T_{\mathbb R})_p(\pi^{-1}(\{\pi(p)\}))\ 
| \ \|\dot{p}\|_p \, < \, 1 \ \}.$$
$\Psi$ is a homeomorphism from $U$ 
to $\cup_{\, p \in \, u(T)} \ B_p(T_{\mathbb R})$. 
For any $t\in T$, 
$\Psi |_t \, 
:= \, \Psi_{\upharpoonright{U_t}}$ 
is a biholomorphic map from $U_t$ to $B_{u(t)}(T_{\mathbb R})$ 
and $\Psi(u(t))=0_{u(t)}$ holds. 

(3) \ \ \  
$$\sup_{\, t\in T, \, \dot{p}\in B_{u(t)}(T_{\mathbb R}) } \ \| \, 
(D(\Psi |_t^{-1}))_{\dot{p}}
 \, \|_{L((T_{\mathbb R})_{u(t)}(\pi^{-1}(\{t\}))
;(T_{\mathbb R})_{\Psi |_t^{-1}(\dot{p})}
(\pi^{-1}(\{t\})))} \ <\ +\infty$$
holds. 

\newpage 

(4) \ \ \ 
For any $r \in [0,1)$ and $\varepsilon>0$, 
there exists $\delta>0$ 
such that for any $t\in T$, $p\in U_t$ and $q\in \pi^{-1}(\{t\})$, 
$$\|\Psi(p)\|_{u(t)} \leq r, \ \ \ d_t(q,p)<\delta$$
$$\Longrightarrow \ \ \ \ \ \ q\in U_t, \ \ \ 
\|\Psi(q)-\Psi(p)\|_{u(t)} < \varepsilon$$
holds. 
\hfill ---

\noindent 
{\bf Definition 3.4} (Coordinate neighborhood of the set of all continuous sections 
defined by a holomorphic normal coordinate neighborhood) : 

Let $T \, (\not=\emptyset)$ be a subset of $X$. 
Let $u \, \in \, \Gamma(M|T)$. 
Let $(U, \Psi)$ to be a holomorphic normal coordinate neighborhood of $u$. 
Then, let 
$$\tilde{U} \ := \ 
\{ \ v\in \Gamma (M|T) \ | \ \exists \, r\in[0,1), \, \forall \, t\in T: 
\ v(t)\in U_t, \, \|\Psi(v(t))\|_{u(t)}\leq r \ \}.$$
A map $\tilde{\Psi}\, :\, \tilde{U}\, 
\rightarrow\, \Gamma_b(u^*(T_{\mathbb R}(M)))$ 
is defined as 
$$(\tilde{\Psi}(v))\, (t)\ :=\ \Psi(v(t)) \ \ \ \ \ \ 
(\, v\in \tilde{U},\, t\in T\, ).$$
$(\tilde{U},\tilde{\Psi})$ 
is called the coordinate neighborhood of $\Gamma (M|T)$
defined by $(U, \Psi)$. 
$\mbox{\ \ \ }$ \hfill ---

\vspace*{0.4em}

The next proposition is the technical main result of this note. 
This proof is given in the next section. 

\vspace*{0.2em}

\noindent 
{\bf Proposition 3.5} 

\noindent (Existence of a holomorphic normal coordinate neighborhood) : 

Let $T\, (\not=\emptyset)$ be a normal subset of $X$. 
Let $u \, \in \, \Gamma(M|T)$. Then, 
there exists a holomorphic normal coordinate neighborhood of $u$. 

\noindent 
{\sf Proof} : 
It is given in the next section. 
\hfill 
$\blacksquare$ 

\vspace*{0.4em}

On the other hand, the main result of this note is as follows. 

\vspace*{0.2em}

\noindent 
{\bf Theorem 3.6} (Main result) : 

Let $T \, (\not=\emptyset)$ be a normal subset of $X$. 
Then, the metric space $\Gamma (M|T)$ 
becomes a $C_b(T)$-manifold 
according to all coordinate neighborhoods of $\Gamma (M|T)$ 
defined by holomorphic normal coordinate neighborhoods. 
\hfill ---

\vspace*{0.2em}

For the proof of this theorem, 
Proposition 3.5 and the following lemma are used. 

\vspace*{0.2em}

\noindent 
{\bf Lemma 3.7} : 

Let $\varepsilon>0$. Let $X$ and $Y$ 
be $n$-dimensional complex inner product spaces. 
Let $f$ be a holomorphic map 
from $\{ \, z\in X \, | \, \|z\|_X<\varepsilon \, \}$ 
to $\{ \, w\in Y \, | \, \|w\|_Y<1 \, \}$. Then, 

\newpage \noindent 
$$\|(Df)_0\|_{L(X;Y)} \ \leq \ \frac{4n^2}{\varepsilon}$$
holds and 
$$z\in X, \ \|z\|_X\leq\frac{1}{4}\varepsilon$$
$$\Longrightarrow
\ \ \ \| \, f(z) \, - \, ( f(0) + (Df)_0(z) ) \, \|_Y
\ \leq \ \frac{16n^3}{\varepsilon^2} \|z\|_X^2$$
holds. 

\noindent 
{\sf Proof} : 
We give it in Section 9. 
\hfill 
$\blacksquare$ 

\noindent 
{\sf Proof of Theorem 3.6} : 

$1^\circ$: \ \ \ 
Let $u \, \in \, \Gamma(M|T)$. 
Let $$\tilde{B}_u\ :=\ \{\, \dot{u} \in  \Gamma_b(u^*(T_{\mathbb R}(M))) \, 
| \, \|\dot{u}\|_u<1 \, \}.$$
Let $(U, \Psi)$ be a holomorphic normal coordinate neighborhood of $u$. 
Then, we show 
$$\tilde{\Psi}(\tilde{U}) \ \subset \ \tilde{B}_u.$$

Let $v\in \tilde{U}$. $v:T\rightarrow U$ 
and $\Psi:U\rightarrow u^*(T_{\mathbb R}(M))$ are continuous. 
Also, there exists $r\in[0,1)$ such that 
$\sup_{t\in T}\|\Psi(v(t))\|_{u(t)} \leq r$ holds. 
Thus, $\tilde{\Psi}(v)=\Psi\circ v \, \in \, \tilde{B}_u$ holds.

$2^\circ$: \ \ \ 
Let $u \, \in \, \Gamma(M|T)$. 
Let $(U, \Psi)$ be a holomorphic normal coordinate neighborhood of $u$. 
A map $\tilde{\Psi^{-1}}\, :\, \tilde{B}_u\, 
\rightarrow\, \Gamma(M|T)$ is defined as 
$$(\tilde{\Psi^{-1}}(\dot{u}))\, (t)\ :=\ \Psi^{-1}(\dot{u}(t)) \ \ \ \ \ \ 
(\, \dot{u}\in \tilde{B}_u,\, t\in T\, ).$$ 
Then, we show 
$$ \tilde{\Psi^{-1}}(\tilde{B}_u) \ \subset \ \tilde{U}.$$

Let $\dot{u}\in \tilde{B}_u$. $\dot{u}\, :\ T\ \rightarrow \  
\cup_{\, p \in \, u(T)} \ B_p(T_{\mathbb R})$ and  
$\Psi^{-1}\, :\ \cup_{\, p \in \, u(T)} \ B_p(T_{\mathbb R}) \ 
\ \rightarrow \ U$ are continuous. 
Also, $\|\dot{u}\|_u\in [0,1)$ and $\|\dot{u}(t)\|_{u(t)}\leq \|\dot{u}\|_u$ hold. 
Thus, $\tilde{\Psi^{-1}}(\dot{u})=\Psi^{-1}\circ \dot{u} \, \in \, \tilde{U}$ holds. 

$3^\circ$: \ \ \ 
From $1^\circ$ and $2^\circ$, 
$$\tilde{\Psi^{-1}}\circ\tilde{\Psi}\, =\, 1_{\tilde{U}}, \ \ \ 
\tilde{\Psi}\circ\tilde{\Psi^{-1}}\, =\, 1_{\tilde{B}_u} $$
hold. 

$4^\circ$: \ \ \ 
We show that $\tilde{\Psi^{-1}}$ is Lipschitz continuous. 

Let $\dot{u_0}, \dot{u_1} \, \in \, \tilde{B}_u$. Then, 
for any $s\in [0,1]$, 
$$\dot{u_s}\, :=\, \dot{u_0}+s(\dot{u_1}-\dot{u_0})
\, = \, (1-s)\dot{u_0}+s\dot{u_1} 
\ \in \ \tilde{B}_u$$
holds. Thus, 

\newpage \noindent 
$$d_t((\tilde{\Psi^{-1}}(\dot{u_0}))(t),(\tilde{\Psi^{-1}}(\dot{u_1}))(t))$$
$$=\ d_t(\Psi^{-1}(\dot{u_0}(t)),\Psi^{-1}(\dot{u_1}(t)))$$
$$\leq \ \int_0^1\ \|\, \frac{d}{ds}(\Psi^{-1}(\dot{u_s}(t)))\, \|_{\Psi^{-1}(\dot{u_s}(t))} \ ds$$
$$= \ \int_0^1\ \|\, (D(\Psi |_t^{-1}))_{\dot{u_s}(t)} \, 
(\dot{u_1}(t)-\dot{u_0}(t))
\, \|_{\Psi |_t^{-1}(\dot{u_s}(t))} \ ds$$
$$\leq \ (\, \int_0^1  \, 
\| (D(\Psi |_t^{-1}))_{\dot{u_s}(t)}
\|_{L((T_{\mathbb R})_{u(t)}(\pi^{-1}(\{t\}))
;(T_{\mathbb R})_{\Psi |_t^{-1}
(\dot{u_s}(t))}(\pi^{-1}(\{t\})))} \, ds\, ) \ 
\|\dot{u_1}(t)-\dot{u_0}(t)\|_{u(t)}$$
holds. Therefore, 
$$d_T(\tilde{\Psi^{-1}}(\dot{u_0}),\tilde{\Psi^{-1}}(\dot{u_1}))$$
$$\leq \ 
(\ \sup_{\, t\in T, \, \dot{p}\in B_{u(t)}(T_{\mathbb R}) } \ \|\, 
(D(\Psi |_t^{-1}))_{\dot{p}}
 \, \|_{L((T_{\mathbb R})_{u(t)}(\pi^{-1}(\{t\}))
;(T_{\mathbb R})_{\Psi |_t^{-1}
(\dot{p})}(\pi^{-1}(\{t\})))} \ ) \ 
\|\dot{u_1}-\dot{u_0}\|_{u}$$
holds. 

$5^\circ$: \ \ \ 
We show that $\tilde{\Psi}$ is continuous. 

Let $v\in\tilde{U}$ and $\varepsilon>0$. 
Then, there exists $r\in [0,1)$ such that for any $t\in T$, 
$$v(t) \ \in \ U_t, \ \ \ \|\, \Psi(v(t))\, \|_{u(t)} \ \leq \ r$$
hold. Therefore, there exists $\delta>0$ such that 
for any $t\in T$ and $q\in \pi^{-1}(\{t\})$,  
$$d_t(q,v(t))<\delta$$
$$\Longrightarrow \ \ \ \ \ \ q\in U_t, \ \ \ 
\|\Psi(q)-\Psi(v(t))\|_{u(t)} < \frac{\varepsilon}{2}$$ 
holds. Let $w\in \tilde{U}$ and $d_T(w,v)<\delta$. Then, 
because of $d_t(w(t),v(t))<\delta$, 
$$\| \, \Psi(w(t)) - \Psi(v(t)) \, \|_{u(t)} \ < \ \frac{\varepsilon}{2}$$
holds. So, 
$$\| \, \tilde{\Psi}(w) - \tilde{\Psi}(v) \, \|_{u} \ 
\leq \ \frac{\varepsilon}{2} \ < \ \varepsilon$$
holds.

$6^\circ$: \ \ \ 
We show that $\tilde{U}$ is an open set of $\Gamma(M|T)$. 

Let $v\in\tilde{U}$. Then, there exists $r\in [0,1)$ such that 
for any $t\in T$, 

\newpage \noindent 
$$v(t) \ \in \ U_t, \ \ \ \|\, \Psi(v(t))\, \|_{u(t)} \ \leq \ r$$
hold. Thus, there exists $\varepsilon>0$ such that 
for any $t\in T$ and $q\in \pi^{-1}(\{t\})$, 
$$d_t(q,v(t))<\varepsilon$$
$$\Longrightarrow \ \ \ \ \ \ q\in U_t, \ \ \ 
\|\Psi(q)-\Psi(v(t))\|_{u(t)} < \frac{1-r}{2}$$
holds. Let $w\in \Gamma(M|T)$ and $d_T(w,v)<\varepsilon$. Then, 
because of $d_t(w(t),v(t))<\varepsilon$, 
$$w(t)\in U_t, \ \ \ \| \, \Psi(w(t)) - \Psi(v(t)) \, \|_{u(t)} \ < \ \frac{1-r}{2}$$
hold. Therefore, furthermore, 
$$\| \, \Psi(w(t)) \, \|_{u(t)} $$
$$\leq \ \| \, \Psi(w(t)) - \Psi(v(t)) \, \|_{u(t)} \, 
+ \, \| \, \Psi(v(t)) \, \|_{u(t)}$$
$$ < \ \frac{1-r}{2}+r \ = \ \frac{1+r}{2}$$
holds, while $\frac{1+r}{2}\in[0,1)$ holds. 
Hence, $w\in \tilde{U}$ holds. 

$7^\circ$: \ \ \ $\tilde{B}_u$ is an open set of $\Gamma_b(u^*(T_{\mathbb R}(M)))$. 
On the other hand, from $\Psi(u(t))=0_{u(t)}$, we obtain $u \in \tilde{U}$. 
Therefore, by Proposition 3.5, $\Gamma (M|T)$ 
is a topological manifold according to all coordinate neighborhoods 
defined by holomorphic normal coordinate neighborhoods.

$8^\circ$: \ \ \ 
We show that each coordinate transformation 
is $C_b(T)$-differentiable (i.e., it is Frechet differentiable 
and its Frechet derivatives are $C_b(T)$-linear).

Let $(U,\Psi)$ be a holomorphic normal coordinate neighborhood of $u$. 
Let $(V,\Phi)$ be a holomorphic normal coordinate neighborhood of $v$. 
Let $\dot{u} \in \tilde{\Psi}(\tilde{U}\cap\tilde{V})$. Then, 
there exists $\varepsilon>0$ such that 
for any $\dot{h}\in\Gamma_b(u^*(T_{\mathbb R}(M)))$, 
$$\|\dot{h}\|_u \, < \, \varepsilon 
\ \ \ \Longrightarrow \ \ \ \dot{u}+\dot{h} 
\, \in \, \tilde{\Psi}(\tilde{U}\cap\tilde{V})$$
holds. Also, because $T$ is normal, if $\|\dot{p}\|_{u(t)}<\varepsilon$ holds, then 
there exists $\dot{h}\in\Gamma_b(u^*(T_{\mathbb R}(M)))$ such that 
$$\|\dot{h}\|_u \, < \, \varepsilon,  
\ \ \ \dot{h}(t)=\dot{p}$$
hold. (For this, since we can take a local trivialization coordinate system 
of the vector bundle $u^*(T_{\mathbb R}(M))$ 
and make the Schmidt orthogonalization on the canonical basis at each $t$ 
to obtain an isomorphism to the product bundle of the complex inner product space, 
it is only necessary to paste the constant section through $\dot{p}$ 
and the zero section by a partition of unity.) Hence, 
$$\|\dot{p}\|_{u(t)} \, < \, \varepsilon \ \ \ 
\Longrightarrow \ \ \ 
\dot{u}(t)+\dot{p} \,  
\in \, \Psi |_t(U_t\cap V_t)$$
holds. So, by Lemma 3.7, 
\begin{equation}\|(D(\Phi |_t
\circ \Psi |_t^{-1}))_{\dot{u}(t)}\|
_{L((T_{\mathbb R})_{u(t)}(\pi^{-1}(\{t\}));
(T_{\mathbb R})_{v(t)}(\pi^{-1}(\{t\})))}
\leq \frac{4n^2}{\varepsilon}\end{equation}
holds and also, 
\begin{equation}\|\dot{p}\|_{u(t)} \, \leq \, \frac{1}{4}\varepsilon
\end{equation}
$$\Longrightarrow$$ 
$$\| \, (\Phi |_t
\circ \Psi |_t^{-1})
(\dot{u}(t)+\dot{p}) \, 
- \, ( \, (\Phi |_t
\circ \Psi |_t^{-1})
(\dot{u}(t)) \, 
+ \, (D(\Phi |_t
\circ \Psi |_t^{-1}))_{\dot{u}(t)} (\dot{p}) \, )
\, \|_{v(t)}$$
$$\leq \ \frac{16n^3}{\varepsilon^2} \, \|\dot{p}\|_{u(t)}^2$$
holds. 

Now, for 
$\dot{h} \in \Gamma_b(u^*(T_{\mathbb R}(M)))$ 
and $t\in T$, we define 
$$(\tilde{A}(\dot{h}))(t) 
\ := \ (D(\Phi |_t
\circ \Psi |_t^{-1}))_{\dot{u}(t)} 
\, (\dot{h}(t)).$$
Then, 
$(\tilde{A}(\dot{h}))(t) 
\, \in \, (T_{\mathbb R})_{v(t)}(\pi^{-1}(\{t\}))$
holds. Furthermore, in virtue of Lemma 1.11, 
the section 
$$t\in T \ \mapsto \ 
(\tilde{A}(\dot{h}))(t) \in v^*(T_{\mathbb R}(M))$$
is continuous. 
(This is because 
taking a respective local trivialization coordinate system 
of the vector bundle $u^*(T_{\mathbb R}(M))$ and $v^*(T_{\mathbb R}(M))$ 
and displaying the map $\dot{p} \, \mapsto \, (\Phi\circ\Psi^{-1})(\dot{p})$ 
in the local coordinates 
to be holomorphic on each fiber and to be continuous, 
so that the map $\dot{p} \, \mapsto \, (D(\Phi |_{\pi(\varpi(\dot{p}))}
\circ \Psi |_{\pi(\varpi(\dot{p}))}^{-1}))_{\dot{p}}$ is also continuous.) 
Also, while it is obvious that $\tilde{A}$ is $C_b(T)$-linear, 
from $(3.1)$, $\tilde{A}$ is a continuous map from $\Gamma_b(u^*(T_{\mathbb R}(M)))$ 
to $\Gamma_b(v^*(T_{\mathbb R}(M)))$. However, 
from $(3.2)$, $\tilde{A}$ is the Frechet derivative 
of $\tilde{\Phi}\circ \tilde{\Psi^{-1}}$ at the point $\dot{u}$. 
\hfill 
$\blacksquare$

\noindent 
{\bf Corollary 3.8} : 

Suppose that $T$ is a paracompact contractible subset of $X$. 
Then, $\Gamma(M|T)$ is an $n$-dimensional $C_b(T)$-manifold. 

\noindent 
{\sf Proof} : 
Let $u\in \Gamma(M|T)$. 
Then, as a continuous complex vector bundle on $T$, 
$u^*(T_{\mathbb R}(M))$ is isomorphic to ${\mathbb C}^n\times T$. 
That is, there exists $\{\dot{u}_k\}_{k=1}^n$ 
such that $\dot{u}_k \, \in \, \Gamma(u^*(T_{\mathbb R}(M)))$ holds 
and for any $t\in T$, $\{\dot{u}_k(t)\}_{k=1}^n$ 
is a basis of $(T_{\mathbb R})_{u(t)}(\pi^{-1}(\{t\}))$. 
We denote the Schmidt orthogonalization 
on $\{\dot{u}_k(t)\}_{k=1}^n$ at each $t\in T$ 
by $\{\dot{e}_k(t)\}_{k=1}^n$. Then, 
$\dot{e}_k \, \in \, \Gamma_b(u^*(T_{\mathbb R}(M)))$ holds 
and for any $t\in T$, $\{\dot{e}_k(t)\}_{k=1}^n$ 
is an orthonormal basis of $(T_{\mathbb R})_{u(t)}(\pi^{-1}(\{t\}))$. 
We define a map $F \, : \, (C_b(T))^n 
\, \rightarrow \, \Gamma_b(u^*(T_{\mathbb R}(M)))$ as 
$$F \, ( \, \{\dot{z}_1(t)\}_{t\in T},\{\dot{z}_2(t)\}_{t\in T},
\cdots,\{\dot{z}_n(t)\}_{t\in T} \, )$$
$$:= \ \{ \, \sum_k \dot{z}_k(t)\dot{e}_k(t) \, \}_{t\in T}.$$ 
Because of $\max_k |z_k| \, \leq \, \|z\|_{{\mathbb C}^n} 
\, \leq \, n \max_k |z_k|$, 
as in Example 2.9, $F$ and $F^{-1}$ 
are continuous $C_b(T)$-linear maps. 
\hfill 
$\blacksquare$ 

\noindent 
{\sf Remark} (Serre-Swan theorem) : 

By making a continuous complex vector bundle $M$ on $X$ 
correspond to the module $\Gamma(M)$ 
of all continuous sections of $M$ on $X$, 
the category of continuous complex vector bundles on $X$ 
is equivalent to the one of finitely generated projective 
$C(X)$-modules ([2, 9, 11, 12]).  
\hfill ---

\noindent 
{\sf Remark} (Manifold on the real commutative algebra $C^\infty(X;{\mathbb R})$) : 

Let $M\rightarrow X$ be a differentiable family of differentiable manifolds. 
The set of all differentiable sections of $M$ on $X$ 
becomes a $C^\infty(X;{\mathbb R})$-manifold ([4, 6, 7, 8]). 
However, $C^\infty(X;{\mathbb R})$ 
is not a real Banach space, but it is a real Frechet space. 
\hfill ---

\noindent 
{\sf Conjecture} (Perhaps not too difficult) : 

For example, similarly, it seems that when $M$ is a compact continuous family 
of almost complex manifolds on $X$, 
an almost complex $C(X;{\mathbb R})$-manifold structure 
is introduced into $\Gamma(M)$. 
\hfill ---

\newpage

\sect{Holomorphic linear connections} 
\ \ \ \ \ Sprays are used 
in order to prove existence of a holomorphic normal coordinate neighborhood. 
However, the global sprays corresponding to Levi-Civita connections may not be holomorphic. 
For each section, holomorphic sprays on its neighborhood 
are constructed as the ones corresponding to convex combinations of trivial connections. 

\vspace*{0.4em}

\noindent 
{\bf Definition 4.1} (Holomorphic linear spray) : 

Let $Z$ be a real $C^\infty$-vector field on the real tangent bundle $T_{\mathbb R}(N)$ 
of an $n$-dimensional complex manifold $N$. Then, 
$Z$ is said to be a holomorphic linear spray on $N$, 
if for any holomorphic local coordinate system $z=(z_1,z_2,\cdots,z_n)$ of $N$, 
there exists a family $\{\Gamma_{i,j}^k\}_{\, i,j,k \, \in \, \{1,2,\cdots,n\}}$ 
of holomorphic functions on the coordinate neighborhood of $z$ such that 
$$\Gamma_{i,j}^k=\Gamma_{j,i}^k \ \ \ \ \ \ (\, i,j,k \, = \, 1, 2, \cdots, n\, )$$
holds and the ordinary differential equation 
that an integral curve 

\noindent 
$\{\, (\dot{z}(s),z(s))\, \}_{\, s\, \in \, (-\varepsilon,+\varepsilon)}$ 
of the vector field $Z$ should satisfy is 
$$\frac{d}{ds}\dot{z}_k 
\ = \ - \, \sum_{i,j} \, \Gamma_{i,j}^k(z_1,z_2,\cdots,z_n) \, \dot{z}_i \, \dot{z}_j,$$
$$\frac{d}{ds}z_k \ = \ \dot{z}_k.$$
\hfill ---

\vspace*{0.4em}

Again, recall that a continuous family $M \, 
(=(M,X,\pi,S,\{\, (\rho_\varphi,M^\prime_\varphi,X^\prime_\varphi)\, \}_{\varphi\in S},r_0))$
with bumps was fixed to be one.

\vspace*{0.2em}

\noindent 
{\bf Definition 4.2} (Trivial spray) : 

(1) \ \ \ 
Let $\varphi\in S$. Let 
$$N_\varphi \ := \ (\varphi,\pi |_\varphi)^{-1}
(D_{\frac{r_0+1}{2}}^n\times \overline{X^\prime_\varphi})$$
$$= \ \varphi^{-1}(D_{\frac{r_0+1}{2}}^n) 
\, \cap \, \pi^{-1}(\overline{X^\prime_\varphi}).$$
$N_\varphi$ is an open set of $\pi^{-1}(\overline{X^\prime_\varphi})$. 
Further, let $t\in X$. Let  
$$N_{\varphi,t} \ := \ N_\varphi \, \cap \, \pi^{-1}(\{t\})$$
$$= \ \left\{\begin{array}{ccr} 
\varphi |_t^{-1}(D_{\frac{r_0+1}{2}}^n) 
& & (t\in \overline{X^\prime_\varphi}), 
\\ \emptyset & & (t\in X\setminus\overline{X^\prime_\varphi}).  
\end{array}\right.$$ 

\newpage \noindent 
$N_{\varphi,t}$ is an open set 
of the $n$-dimensional complex manifold $\pi^{-1}(\{t\})$. 

(2) \ \ \ 
Let $\varphi \in S$ and $t\in X$. 
A real $C^\infty$-vector field $Y_{\varphi,t}$ 
on the real tangent bundle $T_{\mathbb R}(N_{\varphi,t})$ is determined 
as with respect to the holomorphic coordinate system 
$p \, \in \, N_{\varphi,t} 
\ \mapsto \ z=\varphi(p) \, \in \, {\mathbb C}^n$ 
of $N_{\varphi,t}$, the ordinary differential equation 
that an integral curve $(\dot{z},z)$ of $Y_{\varphi,t}$ 
should satisfy is 
$$\frac{d}{ds}\dot{z} \ = \ 0,
\ \ \ \ \ \ \frac{d}{ds}z \ = \ \dot{z}.$$ 
\hfill ---

\noindent 
{\bf Proposition 4.3} : 

(1) \ \ \ The trivial spray $Y_{\varphi,t}$ 
is a holomorphic linear one on $N_{\varphi,t}$. 

(2) \ \ \ Let $\psi \, \in \, S$. 
With respect to the holomorphic local coordinate system 
$p \, \in \, N_{\psi,t}\cap N_{\varphi,t} 
\ \mapsto \ z=\psi(p) \, \in \, {\mathbb C}^n$ of $N_{\varphi,t}$, 
the ordinary differential equation 
that an integral curve $(\dot{z},z)$ of the trivial spray $Y_{\varphi,t}$ 
should satisfy is 
$$\frac{d}{ds}\dot{z}_k 
\ = \ - \, \sum_{i,j} \, (\Gamma(\varphi,\psi))_{i,j}^k(z,t) \, \dot{z}_i \, \dot{z}_j,$$
$$\frac{d}{ds}z_k \ = \ \dot{z}_k,$$ 
provided that ${\varphi |_t}_{(i)}$ 
is the $i$-th component of $\varphi |_t$,  
${\psi |_t}_{(j)}$ is the $j$-th component of $\psi |_t$ and 
$$(\Gamma(\varphi,\psi))_{i,j}^k(z,t)$$
$$:= \ - \ \sum_{l,m} 
\ ((\partial_{z_l}\partial_{z_m}
({\psi |_t}_{(k)}\circ \varphi |_t^{-1}))
((\varphi |_t\circ \psi |_t^{-1})(z)))$$
$$((\partial_{z_i}
({\varphi |_t}_{(l)}\circ \psi |_t^{-1}))
(z))
\ ((\partial_{z_j}
({\varphi |_t}_{(m)}\circ \psi |_t^{-1}))
(z))$$
holds. 

\noindent 
{\sf Proof} : Simple calculation. 
\hfill 
$\blacksquare$ 

\noindent 
{\bf Lemma 4.4} : 

Let $(\Gamma(\varphi,\psi))_{i,j}^k$ 
be the same as Proposition 4.3. Then, the function 
$$(z,t) \, \in \, (\psi,\pi |_\psi)(N_\varphi\cap N_\psi) 
\ \mapsto \ (\Gamma(\varphi,\psi))_{i,j}^k(z,t) \, \in \, {\mathbb C}$$
is continuous. 

\noindent 
Proof: According to Lemma 1.11. 
\hfill 
$\blacksquare$

\newpage \noindent 
{\bf Definition 4.5} (Weighting space) : 

Let $R \, (\not=\emptyset)$ be a subset of $S$. Let 
$$W_R \ := \ \{ \ c \, \in \, [0,1]^S \ | 
\ \sum_{\varphi\in S}c_\varphi=\sum_{\varphi\in R}c_\varphi=1 \ \}.$$ 
$W_R$ is called the weighting space of $R$ 
and an element of $W_R$ is called a weighting of $R$. 

\noindent 
{\sf Remark} : 

$W_R$ is a compact subset of $[0,1]^S$. 
$R_1\subset R_2$ implies $W_{R_1}\subset W_{R_2}$. 
\hfill ---

\noindent 
{\bf Definition 4.6} (Spray determined by a weighting) : 

(1) \ \ \ 
Let $R \, (\not=\emptyset)$ be a subset of $S$. Let 
$$N^R \ := \ \cap_{\varphi \in R} \, N_\varphi 
\ = \ ( \, \cap_{\varphi \in R} \, \varphi^{-1}(D_{\frac{r_0+1}{2}}^n) \, ) 
\, \cap \, \pi^{-1}( \, \cap_{\varphi \in R} \, \overline{X^\prime_\varphi} \, ).$$ 
$N^R$ is an open set of $\pi^{-1}( \, \cap_{\varphi \in R} 
\, \overline{X^\prime_\varphi} \, )$. 
Further, let $t\in X$. Let 
$$N^R_t \ := \ \cap_{\varphi \in R} \, N_{\varphi,t} 
\ = \ N^R\cap\pi^{-1}(\{t\})$$
$$= \ \left\{\begin{array}{ccr} 
\cap_{\varphi \in R} \ ( \, \varphi |_t^{-1}(D_{\frac{r_0+1}{2}}^n) \, )  
& & (t\in \cap_{\varphi \in R} \, \overline{X^\prime_\varphi}), 
\\ \emptyset & & (t\in X\setminus\cap_{\varphi \in R} \, \overline{X^\prime_\varphi}).  \end{array}\right.$$
$N^R_t$ is an open set 
of the $n$-dimensional complex manifold $\pi^{-1}(\{t\})$. 

(2) \ \ \ 
For $c\in W_S$ and $t\in X$, a holomorphic linear spray 
$Y^{c,t}$ on $N^{ c^{-1} ((0,1]) }_t$ is defined as 
$$Y^{c,t} \ := \ \sum_{\varphi\in S} \, c_\varphi Y_{\varphi,t}.$$ 
For a subset $R \, (\not=\emptyset)$ of $S$, 
$c\in W_R$ and $t\in X$, a holomorphic linear spray $Y_R^{c,t}$ 
on $N^R_t$ is defined as 
$$Y_R^{c,t} \ := \ {Y^{c,t}}_{\upharpoonright N^R_t}
\ = \ \sum_{\varphi\in R} 
\, c_\varphi {Y_{\varphi,t}}_{\upharpoonright N^R_t}.$$

\noindent 
{\sf Remark} : 

When $c\in W_R$ holds, $c^{-1}((0,1])\subset R$ 
and $N^R_t\subset N^{c^{-1}((0,1])}_t$ hold. 
\hfill ---

\noindent 
{\bf Definition 4.7} (Coordinate display of the spray) :  

Let $\psi\in R\subset S$. Let 
$$N^{R,\psi} \ := \ (\psi,\pi |_\psi)(N^R),$$
$$U^{R,\psi} \ := \ {\mathbb C}^n 
\, \times \, N^{R,\psi} \, \times \, W_R.$$ 

\newpage \noindent 
$N^{R,\psi}$ is an open set of $D_{\frac{r_0+1}{2}}^n \, \times \, 
( \, \cap_{\varphi \in R} \, \overline{X^\prime_\varphi} \, )$ 
and $U^{R,\psi}$ is an open set of 
$({\mathbb C}^n \times D_{\frac{r_0+1}{2}}^n) 
\, \times \, ( \, ( \, \cap_{\varphi \in R} \, \overline{X^\prime_\varphi} \, )
\, \times \, W_R \, )$. 
Let $\{(\Gamma(\varphi,\psi))_{i,j}^k\}_{\varphi,\psi,i,j,k}$ 
be the same as Proposition 4.3. 
A map 
$$f^{R,\psi} \ : \ U^{R,\psi} \, \rightarrow \, {\mathbb C}^{2n}$$ 
is defined as 
$$f^{R,\psi}_k(\dot{z},z,t,c) \ := \ 
- \, \sum_{i,j} \, ( \, \sum_{\varphi \in R}
\, c_\varphi \, (\Gamma(\varphi,\psi))_{i,j}^k(z,t) \, ) 
\, \dot{z}_i \, \dot{z}_j \ \ \ \ \ \ (k=1,2,\cdots,n),$$
$$f^{R,\psi}_k(\dot{z},z,t,c) \ := \ \dot{z}_{k-n}
\ \ \ \ \ \ (k=n+1,n+2,\cdots,2n).$$
\hfill ---

\noindent 
{\bf Definition 4.8} (Coordinate display of the exponential map) : 

Let $\psi\in R\subset S$. Then, 
the set of all $(\dot{z},z,t,c)\in U^{R,\psi}
\, (={\mathbb C}^n \times N^{R,\psi} \times W_R)$ 
such that there exists $(\dot{w},w):[0,1]\rightarrow 
{\mathbb C}^n\times D_{\frac{r_0+1}{2}}^n$ such that 
$$\frac{d}{ds}(\dot{w},w)=f^{R,\psi}(\dot{w},w,t,c) 
\ \ \ \ \ \ (s\in [0,1]),$$
$$(\dot{w},w)(0)=(\dot{z},z)$$
hold is denoted by $I^{R,\psi}$. 
Also, a map 
$$(\dot{e}^{R,\psi},e^{R,\psi}) 
\, : \, I^{R,\psi}\rightarrow{\mathbb C}^n\times D_{\frac{r_0+1}{2}}^n$$
is defined by $(\dot{e}^{R,\psi},e^{R,\psi})(\dot{z},z,t,c)=(\dot{w},w)(1)$. 
\hfill ---

\noindent 
{\bf Lemma 4.9} (ODE with a continuous parameter) : 

Let $\Lambda$ be a topological space. 
Let $U$ be an open set of ${\mathbb R}^m\times \Lambda$. 
Suppose that $f\, :U \rightarrow {\mathbb R}^m$ 
is a continuous map. 
Suppose that for any $\lambda \in \Lambda$, the map 
$x \, \mapsto \, f(x,\lambda)$ is differentiable. 
Suppose that the map $D_xf \, : \, U\rightarrow L({\mathbb R}^m;{\mathbb R}^m)$ 
is continuous. Further, let $(x_0,\lambda_0)\in U$. 
Suppose that $u_0 \, :[0,1]\rightarrow{\mathbb R}^m$ 
is the solution of the initial value problem 
$$\frac{d}{ds}u \, = \, f(u,\lambda_0), 
\ \ \ \ \ \ u(0) \, = \, x_0.$$ 
Then, for any $\varepsilon>0$, 
there exists an open set $V$ of $U$ 
such that it satisfies the followings. 

\newpage 

(1) \ \ \ 
$(x_0,\lambda_0)\in V$ holds. 

(2) \ \ \ 
For any $(x,\lambda)\in V$, there exists the solution 
$u \, :[0,1]\rightarrow{\mathbb R}^m$ of the initial value problem 
$$\frac{d}{ds}u \, = \, f(u,\lambda), 
\ \ \ \ \ \ u(0) \, = \, x$$
such that 
$$\sup_{s\, \in\, [0,1]}\, \|\, u(s)-u_0(s)\, \|_{{\mathbb R}^m} \ < \ \varepsilon$$
holds. 

\noindent 
{\sf Proof} : 
We give it in Section 10. 
\hfill 
$\blacksquare$

\noindent 
{\bf Lemma 4.10} (Holomorphic ODE) : 

Let $U$ be an open set of ${\mathbb C}^m$. 
Suppose that $f\, :\, U\rightarrow {\mathbb C}^m$ 
is a holomorphic map. 
Let $V$ be an open set of $U$. Suppose that 
a map $w \, : \, [0,1]\times V\rightarrow U$ 
satisfies 
$$\frac{\partial w}{\partial s}(s,z)=f(w(s,z)), 
\ \ \ \ \ \ w(0,z)=z.$$
Then, for any $s\in[0,1]$, the map 
$$z \, \in \, V \ \mapsto \ w(s,z) \, \in \, U$$
is holomorphic. 

\noindent 
{\sf Proof} : 
We give it in Section 11. 
\hfill 
$\blacksquare$

\noindent 
{\bf Proposition 4.11} : 

Let $\psi\in R\subset S$. Then, 
$I^{R,\psi}$ is an open set of $U^{R,\psi}
\, (={\mathbb C}^n \times N^{R,\psi} \times W_R)$. 
The map $(\dot{e}^{R,\psi},e^{R,\psi})$ is continuous. 
For any $t \, \in \, \cap_{\varphi \in R} \, \overline{X^\prime_\varphi}$  
and $c \, \in \, W_R$, the map 
$(\dot{z},z) \, \mapsto \, (\dot{e}^{R,\psi},e^{R,\psi})(\dot{z},z,t,c)$ 
is holomorphic. For any $(z,t)\in N^{R,\psi}$ and $c\in W_R$, 
$$(0,z,t,c) \in I^{R,\psi},$$
$$(\dot{e}^{R,\psi},e^{R,\psi})(0,z,t,c)=(0,z),$$
$$(D_{\dot{z}} \, \dot{e}^{R,\psi}) \, (0,z,t,c) 
\, = \, 1_{{\mathbb C}^n} \, = \, 1_{{\mathbb R}^{2n}},$$
$$(D_{\dot{z}} \, e^{R,\psi}) \, (0,z,t,c) 
\, = \, 1_{{\mathbb C}^n} \, = \, 1_{{\mathbb R}^{2n}}$$
hold. 

\noindent 
{\sf Proof} : 
In virtue of Lemma 4.4 and Lemma 4.9, 
$I^{R,\psi}$ is an open set of $U^{R,\psi}$ 
and $(\dot{e}^{R,\psi},e^{R,\psi})$ is continuous. 
In virtue of Lemma 4.10, for any 
$t \, \in \, \cap_{\varphi \in R} \, \overline{X^\prime_\varphi}$ 
and $c \, \in \, W_R$, 
$(\dot{z},z) \, \mapsto \, (\dot{e}^{R,\psi},e^{R,\psi})(\dot{z},z,t,c)$ 
is holomorphic.  
The others follow because 
$(\dot{e}^{R,\psi},e^{R,\psi}) \, (\dot{z},z,t,c)$ 
is the value at the time $1$ of the integral curve of the splay 
with its value $(\dot{z},z)$ at the time $0$. 
\hfill 
$\blacksquare$

\noindent 
{\bf Lemma 4.12} (Inverse function theorem with a continuous parameter) : 

Let $\Lambda$ be a topological space. 
Let $U$ be an open set of ${\mathbb R}^m\times \Lambda$. 
Suppose that $f \, : \, U \rightarrow {\mathbb R}^m$ 
is a continuous map. Suppose that 
for any $\lambda \in \Lambda$, the map 
$x \, \mapsto \, f(x,\lambda)$ is differentiable. 
Suppose that the map $D_xf \, : \, U\rightarrow L({\mathbb R}^m;{\mathbb R}^m)$ 
is continuous. Further, suppose that 
$\Theta$ is a compact subset of $\Lambda$. 
Suppose that $\{0\}\times\Theta \, \subset \, U$ and 
$$\lambda \, \in \, \Theta \ \ \ \Longrightarrow 
\ \ \ (D_xf) \, (0,\lambda) \ = \ 1_{ \, {\mathbb R}^m    }$$
hold. Then, for any $\varepsilon>0$, 
there exist $\delta>0$, an open set $\Lambda^\prime$ of $\Lambda$, 
an open set $U^\prime$ of $U$ 
and an open set $V^\prime$ of ${\mathbb R}^m\times \Lambda$ 
such that they satisfy the followings. 

(1) \ \ \ 
$$\Theta \ \subset \ \Lambda^\prime,$$ 
$$\{ \ x \in {\mathbb R}^m \ | \  \|x\|_{{\mathbb R}^m}<\delta \ \}
\ \times \ \Lambda^\prime 
\ \ \ \subset \ \ \ U^\prime 
\ \ \ \subset \ \ \ {\mathbb R}^m \ \times \ \Lambda^\prime,$$ 
$$\cup_{ \, \lambda \, \in \, \Lambda^\prime} \ 
( \, \{  \, y \in {\mathbb R}^m \, | \, \|y-f(0,\lambda)\|_{{\mathbb R}^m}<\delta \, \} 
\, \times \, \{\lambda\}  \, ) \ \ \ \subset \ \ \ V^\prime 
\ \ \ \subset \ \ \ {\mathbb R}^m \ \times \ \Lambda^\prime,$$
$$\sup_{(x,\lambda) \, \in \, U^\prime}
\ \| \, (D_xf) \, (x,\lambda) 
\, - \, 1_{ \, {\mathbb R}^m } \, \|_{L({\mathbb R}^m;{\mathbb R}^m)} 
\ \ \ \leq \ \ \ \varepsilon$$
hold. Homeomorphically, the map $(x,\lambda) \, \in \, U \ \mapsto 
\ (f(x,\lambda),\lambda) \, \in \, {\mathbb R}^m\times \Lambda$ 
maps $U^\prime$ to $V^\prime$. 

(2) \ \ \ There exists a continuous map $g \, : \, V^\prime\rightarrow{\mathbb R}^m$ 
such that it satisfies the followings. For any $(y,\lambda) \, \in \, V^\prime$,  
$$( \, g(y,\lambda), \, \lambda \, ) \ \in \ U^\prime, 
\ \ \ f( \, g(y,\lambda), \, \lambda \, ) \ = \ y$$
hold. For any $\lambda \in \Lambda^\prime$, the map 
$y \, \mapsto \, g(y,\lambda)$ is differentiable. 
The map $D_yg:\ V^\prime\rightarrow L({\mathbb R}^m;{\mathbb R}^m)$ 
is continuous. 
$$\sup_{(y,\lambda) \, \in \, V^\prime}
\ \| \, (D_yg) \, (y,\lambda) \, 
- \, 1_{ \, {\mathbb R}^m } \, \|_{L({\mathbb R}^m;{\mathbb R}^m)} 
\ \ \ \leq \ \ \ \varepsilon$$
holds. 

\noindent 
{\sf Proof} : 
We give it in Section 12. 
\hfill 
$\blacksquare$

\noindent 
{\bf Definition 4.13} : 

(1) \ \ \ 
Let $R \, (\not=\emptyset)$ be a subset of $S$. 
Let $$K^R \ := \ \cap_{\varphi \in R} 
\, (\varphi,\pi |_\varphi)^{-1}
(\overline{D_{r_0}^n}\times \overline{X^\prime_\varphi}).$$ 
$K^R$ is a compact subset of $N^R$. 
Further, let $\psi\in R$ and 
$$K^{R,\psi} \ := \ (\psi,\pi |_\psi)(K^R).$$ 
$K^{R,\psi}$ is a compact subset of $N^{R,\psi}$. 

(2) \ \ \ Let 
$$\mathcal{P}^\prime(S) 
\ := \ \{ \, (R,\psi) 
\, | \, \psi \in R \subset S \, \}.$$ 
\hfill ---

\noindent 
{\bf Lemma 4.14} : 

There exist $\alpha_0>0$, $\beta_0>0$, $\delta_0^a>0$, $\delta_0^b>0$ and 
$\{(U^\prime_{R,\psi},V^\prime_{R,\psi},h^{R,\psi})\}_{(R,\psi)\in\mathcal{P}^\prime(S)}$ 
such that for any $(R,\psi)\in\mathcal{P}^\prime(S)$, they satisfy the followings. 

(1) \ \ \ 
$U^\prime_{R,\psi} \, \subset \, I^{R,\psi}$ holds. 
$U^\prime_{R,\psi}$ is an open set of 
${\mathbb C}^n\times K^{R,\psi}\times W_R$. 
$V^\prime_{R,\psi}$ is an open set of 
$D_{\frac{r_0+1}{2}}^n\times K^{R,\psi}\times W_R$. 
$$\{ \ \dot{z} \in {\mathbb C}^n \ | \  \|\dot{z}\|_{{\mathbb C}^n}
<\delta_0^a \ \}
\ \times \  K^{R,\psi}\times W_R
\ \ \ \subset \ \ \ U^\prime_{R,\psi},$$ 
$$( \, \cup_{ \, (z,t) \, \in \, K^{R,\psi}} \ 
( \, \{  \, w \in {\mathbb C}^n \, | \, \|w-z\|_{{\mathbb C}^n}
<\delta_0^b \, \} 
\, \times \, \{(z,t)\}  \, ) \, ) \, \times \, W_R \ \ \ \subset \ \ \ V^\prime_{R,\psi},$$
$$\sup_{(\dot{z},z,t,c) \, \in \, U^\prime_{R,\psi}}
\ \| \, (D_{\dot{z}} e^{R,\psi}) \, (\dot{z},z,t,c) 
\, \|_{L({\mathbb C}^n;{\mathbb C}^n)} \ \leq \ \alpha_0$$
hold. 
Homeomorphically, the map $(\dot{z},z,t,c) \, \in \, I^{R,\psi} \ \mapsto 
\ (e^{R,\psi}(\dot{z},z,t,c),z,t,c) \, \in \, 
D_{\frac{r_0+1}{2}}^n\times N^{R,\psi}\times W_R$ maps 
$U^\prime_{R,\psi}$ to $V^\prime_{R,\psi}$. 

(2) \ \ \ 
$h^{R,\psi}$ is a continuous map from $V^\prime_{R,\psi}$ 
to ${\mathbb C}^n$. 
For any $(w,z,t,c) \, \in \, V^\prime_{R,\psi}$, 
$$(h^{R,\psi}(w,z,t,c),z,t,c) \ \in \ U^\prime_{R,\psi},$$
$$e^{R,\psi}(h^{R,\psi}(w,z,t,c),z,t,c) \ = \ w$$
hold. For any $(z,t) \in K^{R,\psi}$ and $c\in W_R$, the map 
$w \, \mapsto \, h^{R,\psi}(w,z,t,c)$ is holomorphic. 
$$\sup_{(w,z,t,c) \, \in \, V^\prime_{R,\psi}}
\ \| \, (D_w h^{R,\psi}) \, (w,z,t,c) 
\, \|_{L({\mathbb C}^n;{\mathbb C}^n)} \ \leq \ \beta_0$$
holds. 

\noindent 
{\sf Proof} : 
Since ${\mathcal P}^\prime(S)$ is a finite set, 
it is sufficient to discuss it for each 
$(R,\psi)\in{\mathcal P}^\prime(S)$. 
Let $(R,\psi)\in{\mathcal P}^\prime(S)$. As we set 

\newpage \noindent 
$$\Theta \ := \ \Lambda \ := \ K^{R,\psi}\times W_R,$$
$$U \ := \ I^{R,\psi}\cap({\mathbb C}^n\times K^{R,\psi}\times W_R)
\ = \ I^{R,\psi}\cap({\mathbb R}^{2n}\times \Lambda),$$
$$f \ := \ {e^{R,\psi}}_{\upharpoonright U},$$ 
in virtue of Lemma 1.11 and Proposition 4.11, we can apply Lemma 4.12. 

$\mbox{\ \ \ }$  \hfill 
$\blacksquare$

\vspace*{0.2em}

Hereafter, we fix $\alpha_0>0$, $\beta_0>0$, $\delta_0^a>0$, $\delta_0^b>0$ and 
$\{(U^\prime_{R,\psi},V^\prime_{R,\psi},h^{R,\psi})\}_{(R,\psi)\in\mathcal{P}^\prime(S)}$
given by the above lemma to be ones. 

\vspace*{0.2em}

\noindent 
{\bf Definition 4.15} : 

Let $\{(U^\prime_{R,\psi},V^\prime_{R,\psi},h^{R,\psi})\}_{(R,\psi)\in\mathcal{P}^\prime(S)}$ 
be the same as Lemma 4.14. 
Let $\psi \, \in \, R \, \subset \, S$. 
Further, suppose that $T^\prime$ is a subset of $\overline{X^\prime_\psi}$, 
$\chi^\prime$ is a continuous map from $T^\prime$ to $W_R$ 
and $v^\prime$ is a continuous map from $T^\prime$ to $\overline{D_{r_0}^n}$ 
such that for any $t\in T^\prime$, 
$(v^\prime(t),t) \, \in \, K^{R,\psi}$ 
holds. Then, let 
$$U^\prime_{R,\psi,v^\prime,\chi^\prime}
\ := \ \{ \, (\dot{z},t) \in {\mathbb C}^n \times T^\prime 
\, | \, (\dot{z},v^\prime(t), t, \chi^\prime(t)) 
\in U^\prime_{R,\psi} \, \},$$
$$V^\prime_{R,\psi,v^\prime,\chi^\prime}
\ := \ \{ \, (w,t) \in D_{\frac{r_0+1}{2}}^n \times T^\prime 
\, | \, (w,v^\prime(t), t, \chi^\prime(t)) 
\in V^\prime_{R,\psi} \, \}.$$
Further, a map $e^{R,\psi,v^\prime,\chi^\prime}$ 
from $U^\prime_{R,\psi,v^\prime,\chi^\prime}$
to $D_{\frac{r_0+1}{2}}^n$ is defined as 
$$e^{R,\psi,v^\prime,\chi^\prime}(\dot{z},t) 
\ := \ e^{R,\psi}(\dot{z},v^\prime(t), t, \chi^\prime(t))$$ 
and a map $h^{R,\psi,v^\prime,\chi^\prime}$ 
from $V^\prime_{R,\psi,v^\prime,\chi^\prime}$ 
to ${\mathbb C}^n$ is defined as 
$$h^{R,\psi,v^\prime,\chi^\prime}(w,t) 
\ := \ h^{R,\psi}(w,v^\prime(t), t, \chi^\prime(t)).$$
\hfill ---

\noindent 
{\bf Lemma 4.16} : 

Let $\alpha_0>0$, $\beta_0>0$, $\delta_0^a>0$ and $\delta_0^b>0$ 
be the same as Lemma 4.14. Let $\psi \, \in \, R \, \subset \, S$. 
Further, suppose that $T^\prime$ is a subset of $\overline{X^\prime_\psi}$,  
$\chi^\prime$ is a continuous function from $T^\prime$ to $W_R$ 
and $v^\prime$ is a continuous function from $T^\prime$ to $\overline{D_{r_0}^n}$ 
such that for any $t\in T^\prime$, 
$(v^\prime(t),t) \, \in \, K^{R,\psi}$
holds. Then, the followings hold.  

$U^\prime_{R,\psi,v^\prime,\chi^\prime}$ 
is an open set of ${\mathbb C}^n\times T^\prime$. 
$V^\prime_{R,\psi,v^\prime,\chi^\prime}$ 
is an open set of $D_{\frac{r_0+1}{2}}^n\times T^\prime$. 
$$\{ \ \dot{z} \in {\mathbb C}^n \ | \  \|\dot{z}\|_{{\mathbb C}^n}
<\delta_0^a \ \}
\ \times \ T^\prime \ \ \ \subset \ \ \ 
U^\prime_{R,\psi,v^\prime,\chi^\prime},$$ 
$$\cup_{ \, t \, \in \, T^\prime} \ 
( \, \{  \, w \in {\mathbb C}^n \, | \, \|w-v^\prime(t)\|_{{\mathbb C}^n}
<\delta_0^b \, \} 
\, \times \, \{t\}  \, )
\ \ \ \subset \ \ \ V^\prime_{R,\psi,v^\prime,\chi^\prime}$$
hold. For any $t \in T^\prime$, the maps 
$\dot{z} \, \mapsto \, e^{R,\psi,v^\prime,\chi^\prime}(\dot{z},t)$ 
and $w \, \mapsto \, h^{R,\psi,v^\prime,\chi^\prime}(w,t)$ 
are holomorphic. 
$$\sup_{(\dot{z},t) \, \in \, U^\prime_{R,\psi,v^\prime,\chi^\prime}}
\ \| \, (D_{\dot{z}} e^{R,\psi,v^\prime,\chi^\prime}) \, (\dot{z},t) 
\, \|_{L({\mathbb C}^n;{\mathbb C}^n)} \ \leq \ \alpha_0,$$
$$\sup_{(w,t) \, \in \, V^\prime_{R,\psi,v^\prime,\chi^\prime}}
\ \| \, (D_w h^{R,\psi,v^\prime,\chi^\prime}) \, (w,t) 
\, \|_{L({\mathbb C}^n;{\mathbb C}^n)} \ \leq \ \beta_0$$
hold. Homeomorphically, 
the map $(\dot{z},t) \ \mapsto 
\ (e^{R,\psi,v^\prime,\chi^\prime}(\dot{z},t),t)$ 
maps $U^\prime_{R,\psi,v^\prime,\chi^\prime}$ 
to $V^\prime_{R,\psi,v^\prime,\chi^\prime}$. 
For any $(w,t) \, \in \, V^\prime_{R,\psi,v^\prime,\chi^\prime}$, 
$$(h^{R,\psi,v^\prime,\chi^\prime}(w,t),t) \ \in \ U^\prime_{R,\psi,v^\prime,\chi^\prime},$$ 
$$e^{R,\psi,v^\prime,\chi^\prime}(h^{R,\psi,v^\prime,\chi^\prime}(w,t),t) \ = \ w$$
hold. 

\noindent 
{\sf Proof} : It follows from Lemma 4.14.
\hfill 
$\blacksquare$

\noindent 
{\bf Lemma 4.17} : 

There exist $C_1^a>0$ and $C_1^b>0$ 
such that the following holds.  
For any $\varphi \in S$, 
$p\in (\varphi,\pi |_\varphi)^{-1}
(\overline{D_{\frac{r_0+1}{2}}^n} \times \overline{X^\prime_\varphi})$ 
and $\dot{p}\in (T_{\mathbb R})_p(\pi^{-1}(\{\pi(p)\}))$, 
$$\|\dot{p}\|_{\varphi,p}
\ \leq \ C_1^a \, \|\dot{p}\|_p,
\ \ \ \ \ \ \|\dot{p}\|_p
\ \leq \ C_1^b \, \|\dot{p}\|_{\varphi,p}$$
hold. 

\noindent 
{\sf Proof} : 
Since $S$ is a finite set, it follows from Lemma 1.14. 
\hfill 
$\blacksquare$

\noindent 
{\bf Lemma 4.18} : 

There exist $\delta_2>0$ and $C_2>0$ such that the following holds. 
For any $\varphi \in S$, 
$p\in (\varphi,\pi |_\varphi)^{-1}
(\overline{D_{\frac{r_0+1}{2}}^n} \times \overline{X^\prime_\varphi})$ 
and $q \in \pi^{-1}(\{\pi(p)\})$, 
$$d_{\pi(p)}(q,p) \, < \, \delta_2$$
$$\Longrightarrow$$ 
$$q \, \in \, M_\varphi, 
\ \ \ \|\varphi(q)-\varphi(p)\|_{{\mathbb C}^n} 
\, \leq \, C_2 \, d_{\pi(p)}(q,p)$$
holds. 

\noindent 
{\sf Proof} : 
Since $S$ is a finite set, it follows from Lemma 1.16. 
\hfill 
$\blacksquare$

\noindent 
{\sf Proof of Proposition 3.5} :

$0^\circ$: \ \ \ Let 
$\alpha_0, \beta_0, \delta_0^a, \, \delta_0^b, \, \delta_2, 
\, C_1^a, \, C_1^b, \, C_2$ and 
$\{(U^\prime_{R,\psi},V^\prime_{R,\psi},h^{R,\psi})\}_{(R,\psi)\in\mathcal{P}^\prime(S)}$ 
be the same as Lemma 4.14, Lemma 4.17 and Lemma 4.18. 

In order to define a holomorphic linear spray for each $t\in T$, 
we first select a $t$-dependent weighting (i.e., a partition of unity) $\chi$ as follows. 
Because $\{M^\prime_\varphi\}_{\varphi \in S}$ is a finite open covering of $M$, 
$\{u^{-1}(M^\prime_\varphi)\}_{\varphi \in S}$ is a finite open covering of $T$.  
Because $T$ is normal, there exists $\chi 
\, (:=\{\chi_\varphi(t)\}_{ \, \varphi \in S, \, t\in T})$ such that 
$\chi_\varphi$ is a non-negative continuous function on $T$ and 
$${\rm supp}(\chi_\varphi) \ \subset \ u^{-1}(M^\prime_\varphi), 
\ \ \ \ \ \ \sum_{\varphi \in S} \, \chi_\varphi \ = \ 1$$
hold. 

$1^\circ$: \ \ \ 
For $t\in T$, we define 
$$R_t \ := \ \{ \, \varphi\in S 
\, | \, t\in{\rm supp}(\chi_\varphi) \, \}.$$
Then, 
$$\emptyset \ \not= \ R_t \ \subset \ S,$$
$$\chi(t) \, (:=\{\chi_\varphi(t)\}_{\varphi\in S}) \ \in \ W_{R_t},$$ 
$$u(t) \, \in \, \cap_{\varphi \in R_t} \, M^\prime_\varphi 
\, = \, \cap_{\varphi \in R_t} \, (\varphi,\pi |_\varphi)^{-1}
(D_{r_0}^n\times X^\prime_\varphi) 
\ \subset \ K^{R_t}$$
hold.

$2^\circ$: \ \ \ 
For $p\in M$, we define  
$$B^\prime_p(T_{\mathbb R}) 
\ : \, = \ \{\ \dot{p} \, \in \, (T_{\mathbb R})_p(\pi^{-1}(\{\pi(p)\}))\ 
| \ \|\dot{p}\|_p \, < 
\, \frac{\min \{\delta_0^a,\frac{\delta_0^b}{\alpha_0}\}}{C_1^a} \ \}.$$
For $\dot{p} \, \in \, \cup_{p\in u(T)} \, B^\prime_p(T_{\mathbb R})$, 
we define ${\rm exp}(\dot{p}) \, \in \, \pi^{-1}(T)$ as follows.

Let $\dot{p} \, \in \, \cup_{p\in u(T)} \, B^\prime_p(T_{\mathbb R})$. 
Then, there uniquely exists $t\in T$ such that $\dot{p} \, \in 
\, B^\prime_{u(t)}(T_{\mathbb R})$ holds. In virtue of $1^\circ$, 
there exists $(R,\psi)\in {\mathcal P}^\prime(S)$ such that 
$$\chi(t) \, \in \, W_R, \ \ \ u(t) \, \in \, K^R$$
hold. Because of $K^R \, \subset 
\, (\psi,\pi |_\psi)^{-1}
(\overline{D_{r_0}^n}\times\overline{X^\prime_\psi})$, 
by Lemma 4.17, 
$$\|(D(\psi |_t))_{u(t)}(\dot{p})\|_{{\mathbb C}^n}
\, \leq \, C_1^a\|\dot{p}\|_{u(t)} 
\, < \, \min \{\delta_0^a,\frac{\delta_0^b}{\alpha_0}\}$$
holds. Therefore, from Lemma 4.14, 

\newpage \noindent 
$$( \, (D(\psi |_t))_{u(t)}(\dot{p}), \, 
\psi(u(t)), \, t, \, \chi(t) \, ) \, \in \, U^\prime_{R,\psi}$$
holds. Hence, $\dot{p}$ belongs to the domain of the exponential map 
for the spray $Y^{\chi(t),t}_R$ on $N^R_t$. 
So, further, it does to the one for $Y^{\chi(t),t}$. 
Now, we define ${\rm exp}(\dot{p})$ 
as the value of the exponential map for the spray $Y^{\chi(t),t}$ 
at $\dot{p}$. Then, 
$${\rm exp}(\dot{p}) \ \in \ \pi^{-1}(\{t\}) 
\ \subset \ \pi^{-1}(T)$$
holds. 

$3^\circ$: \ \ \ We define 
$$U \ := \ {\rm exp}(\cup_{p\in u(T)} \, B^\prime_p(T_{\mathbb R})).$$
$U$ is a subset of $\pi^{-1}(T)$. ${\rm exp}$ is a map 
whose domain is $\cup_{p\in u(T)} \, B^\prime_p(T_{\mathbb R})$ 
and whose range is $U$. 

$4^\circ$: \ \ \ 
For $t_0\in T$, we define 
$$T_{t_0} \ := \ ( \, \cap_{\psi\in R_{t_0}} \, (u^{-1}(M^\prime_\psi)) \, )
\ \cap \ ( \, \cap_{\psi \in S\setminus R_{t_0}} 
\, (T\setminus {\rm supp}(\chi_\psi)) \, ).$$ 
$T_{t_0}$ is an open set of $T$. $t_0\in T_{t_0}$ holds. 
Thus, in particular, $\{T_{t_0}\}_{t_0\in T}$ 
is an open covering of $T$. Also, 
$$t \, \in \, T_{t_0}$$
$$\Longrightarrow \ \ \ \ \ \ \chi(t) \ \in \ W_{R_{t_0}}, 
\ \ \ u(t) \, \in \, \cap_{\varphi \in R_{t_0}} \, M^\prime_\varphi 
\ \subset \ K^{R_{t_0}}$$
holds. Therefore, if $\psi \, \in \, R_{t_0}$ holds, then 
$T_{t_0}$ is a subset of $\overline{X^\prime_\psi}$,  
$\chi_{\upharpoonright T_{t_0}}$ is a continuous map from $T_{t_0}$ to $W_R$ 
and $\psi\circ (u_{\upharpoonright T_{t_0}})$ is a continuous map 
from $T_{t_0}$ to $\overline{D_{r_0}^n}$ such that 
for any $t\in T_{t_0}$, 
$((\psi\circ (u_{\upharpoonright T_{t_0}}))(t),t) \, \in \, K^{R_{t_0},\psi}$
holds. 

$5^\circ$: \ \ \ 
For a subset $T^\prime$ of $T$, we define 
$$U^{T^\prime} \ := \ U\cap \pi^{-1}(T^\prime)
\ = \ {\rm exp}(\cup_{p\in u(T^\prime)} \, B^\prime_p(T_{\mathbb R})),$$
$${\rm exp} |_{T^\prime}
\ := \ {\rm exp}_{\upharpoonright \, 
\cup_{p\in u(T^\prime)} \, B^\prime_p(T_{\mathbb R})}.$$
Let $t_0\in T$. 
We examine the map ${\rm exp} |_{T_{t_0}}$ 
and its range $U^{T_{t_0}}$ 
and the map ${\rm exp} |_{\{t_0\}}$ 
and its range $U^{\{t_0\}}$.

By $1^\circ$, there exists $\psi \, \in \, R_{t_0}$. 
In virtue of $4^\circ$, 
$u(T_{t_0}) \, \subset \, \cap_{\varphi \in R_{t_0}}M^\prime_\varphi 
\, \subset \, M_\psi$ holds. Now, 
because $\cup_{p\in u(T_{t_0})} \, B^\prime_p(T_{\mathbb R})$ 
is an open set of 
$u^*(T_{\mathbb R}(M))|T_{t_0} \, (\, = \, T_{\mathbb R}(M)|u(T_{t_0})
\, = \, \cup_{p\in u(T_{t_0})} (T_{\mathbb R})_p(\pi^{-1}(\{\pi(p)\})) \, )$, 
the trivialization coordinate system 
$$\dot{p} \, \in \, u^*(T_{\mathbb R}(M))|T_{t_0}
\ \mapsto \ 
(  \,  (D(\psi |_{\pi(\varpi(\dot{p}))}))_{\varpi(\dot{p})}(\dot{p})    , 
\, \pi(\varpi(\dot{p}))  \, ) \, \in \, {\mathbb C}^n\times T_{t_0}$$
of the vector bundle $u^*(T_{\mathbb R}(M))|T_{t_0}$ on $T_{t_0}$
maps $\cup_{p\in u(T_{t_0})} \, B^\prime_p(T_{\mathbb R})$ 
to an open set of ${\mathbb C}^n\times T_{t_0}$. That is, 
as we denote the image of the set $\cup_{p\in u(T_{t_0})} \, B^\prime_p(T_{\mathbb R})$ 
by $G_{t_0}$, $G_{t_0}$ is an open set of ${\mathbb C}^n\times T_{t_0}$. 
We show $G_{t_0} \, \subset \, U^\prime_{R_{t_0},\psi,
\chi_{\upharpoonright T_{t_0}},\psi\circ(u_{\upharpoonright T_{t_0}})}$. 
As $t\in T_{t_0}$ and $\dot{p}\in B^\prime_{u(t)}(T_{\mathbb R})$ hold, 
because of $u(t) \, \in \, K^{R_{t_0}} \, \subset 
\, (\psi,\pi |_\psi)^{-1}
(\overline{D_{r_0}^n}\times\overline{X^\prime_\psi})$, 
by Lemma 4.17, 
$$\|(D(\psi |_t))_{u(t)}(\dot{p})\|_{{\mathbb C}^n}
\, \leq \, C_1^a\|\dot{p}\|_{u(t)} 
\, < \, \min \{\delta_0^a,\frac{\delta_0^b}{\alpha_0}\}$$
holds. So, in virtue of $4^\circ$ and Lemma 4.16, 
$$( \, (D(\psi |_t))_{u(t)}(\dot{p}), 
\, t \, ) \, \in \, U^\prime_{R_{t_0},\psi,
\chi_{\upharpoonright T_{t_0}},
\psi\circ(u_{\upharpoonright T_{t_0}})}$$
holds. $G_{t_0} \, \subset \, U^\prime_{R_{t_0},\psi,
\chi_{\upharpoonright T_{t_0}},\psi\circ(u_{\upharpoonright T_{t_0}})}$
holds. From the above, $G_{t_0}$ is an open set of $U^\prime_{R_{t_0},\psi,
\chi_{\upharpoonright T_{t_0}},\psi\circ(u_{\upharpoonright T_{t_0}})}$. 
Hence, as we set 
$$H_{t_0} \ := \ \{ \, ( \, e^{R_{t_0},\psi,
\chi_{\upharpoonright T_{t_0}},\psi\circ(u_{\upharpoonright T_{t_0}})}
(\dot{z},t), \, t \, ) \, \in \, D_{\frac{r_0+1}{2}}^n \times T_{t_0}
\, | \, (\dot{z},t) \, \in \, G_{t_0} \, \},$$
by Lemma 4.16, $H_{t_0}$ is an open set of 
$D_{\frac{r_0+1}{2}}^n \times T_{t_0}$ and homeomorphically  
the map 
$$(\dot{z},t) \ \mapsto \ 
( \, e^{R_{t_0},\psi,
\chi_{\upharpoonright T_{t_0}},\psi\circ(u_{\upharpoonright T_{t_0}})}
(\dot{z},t), \, t \, )$$
maps $G_{t_0}$ to $H_{t_0}$. Therefore, $U^{T_{t_0}}$ 
is an open set of $\pi^{-1}(T_{t_0})$ 
and homeomorphically  the map ${\rm exp} |_{T_{t_0}}$ 
maps to its domain $\cup_{p\in u(T_{t_0})} \, B^\prime_p(T_{\mathbb R})$ 
to its range $U^{T_{t_0}}$. Further, ${\rm exp} |_{\{t_0\}}$ 
is a biholomorphic map from its domain $B^\prime_{u(t_0)}(T_{\mathbb R})$ 
to its range $U^{\{t_0\}}$. Also, 
$$u(t_0) \, = \, {\rm exp} |_{\{t_0\}}(0_{u(t_0)})
\, \in \, U^{\{t_0\}}$$
holds. Also, 
because for $\dot{p}\in B^\prime_{u(t_0)}(T_{\mathbb R})$, 
$$(D({\rm exp} |_{\{t_0\}}))_{\dot{p}}$$
$$= \ (D(\psi |_{t_0}^{-1}))
_{(e^{R_{t_0},\psi,
\chi_{\upharpoonright T_{t_0}},\psi\circ(u_{\upharpoonright T_{t_0}})}
((D(\psi |_{t_0}))_{u(t_0)}(\dot{p}),t_0),t_0)}$$
$$\cdot \, (D_{\dot{z}}e^{R_{t_0},\psi,
\chi_{\upharpoonright T_{t_0}},\psi\circ(u_{\upharpoonright T_{t_0}})})
((D(\psi |_{t_0}))_{u(t_0)}(\dot{p}),t_0)
\, \cdot \, (D(\psi |_{t_0}))_{u(t_0)}$$
holds, from Lemma 4.16 and Lemma 4.17, 

\newpage \noindent 
$$\sup_{\dot{p}\in B^\prime_{u(t_0)}(T_{\mathbb R})} \ 
\| \, (D({\rm exp} |_{\{t_0\}}))_{\dot{p}} \, \|
_{ \, L \, ( \, (T_{\mathbb R})_{u(t_0)}(\pi^{-1}(\{t_0\}))  
\, ; \, (T_{\mathbb R})_{{\rm exp} |_{\{t_0\}}
(\dot{p})}(\pi^{-1}(\{t_0\})) \, )}$$
$$\leq \ C_1^b \, \cdot \, \alpha_0 \, \cdot \, C_1^a$$
holds.

$6^\circ$: \ \ \ 
For any $t_0\in T$, $T_{t_0}$ is an open set of $T$ 
and $t_0\in T_{t_0}$ holds. 
Hence, from $5^\circ$, 
$U$ is an open set of $\pi^{-1}(T)$ 
and homeomorphically the map ${\rm exp}$ 
maps its domain $\cup_{p\in u(T)} \, B^\prime_p(T_{\mathbb R})$ 
to its range $U$. Further, for any $t\in T$, 
$$u(t) \ = \ {\rm exp}(0_{u(t)})
\ \in \ U_t \ ( \, := \, U\cap\pi^{-1}(\{t\}) \, )$$
holds and the map ${\rm exp} |_{\{t\}}$ 
is a biholomorphic map from its domain 
$B^\prime_{u(t)}(T_{\mathbb R})$ 
to its range $U_t$. Also, 
$$\sup_{t\in T, \, \dot{p}\in B^\prime_{u(t)}(T_{\mathbb R})} \ 
\| \, (D({\rm exp} |_{\{t\}}))_{\dot{p}} \, \|
_{ \, L \, ( \, (T_{\mathbb R})_{u(t)}(\pi^{-1}(\{t\}))  
\, ; \, (T_{\mathbb R})_{{\rm exp} |_{\{t\}}
(\dot{p})}(\pi^{-1}(\{t\})) \, )}$$
$$\leq \ \alpha_0 C_1^a C_1^b 
\ < \ +\infty$$
holds. 

$7^\circ$: \ \ \ 
For $q \in U$, we define 
$$\Psi(q) \ := \ \frac{C_1^a}{\min \{\delta_0^a,\frac{\delta_0^b}{\alpha_0}\}}
\, {\rm exp}^{-1}(q).$$ 
From $6^\circ$, $(U,\Psi)$ satisfies the conditions (1), (2) and (3) 
of the definition.

$8^\circ$: \ \ \ 
We show that the condition (4) of the definition is satisfied. 

Let $0\leq r<1$ and $\varepsilon>0$.  
We set $$\delta \ := \ \min \ 
\{ \ \delta_2 , 
\ \frac{\delta_0^b}{C_2} \, (1-r), 
\ \frac{ \min \{\delta_0^a,\frac{\delta_0^b}{\alpha_0}\} }{\beta_0C_1^aC_1^bC_2} \, 
(1-r), 
\ \frac{ \min \{\delta_0^a,\frac{\delta_0^b}{\alpha_0}\} }{\beta_0C_1^aC_1^bC_2} \, 
\varepsilon \ \}.$$
$\delta>0$ holds. 

Suppose $t\in T$, $p\in U_t$, $q\in \pi^{-1}(\{t\})$, 
$\|\Psi(p)\|_{u(t)} \leq r$ and $d_t(q,p)<\delta$. 
We show $q \, \in \, U_t$ and  
$\|\Psi(q)-\Psi(p)\|_{u(t)} \, < \, \varepsilon$. 
From $1^\circ$, there exists $\psi$ such that 
$\psi \, \in \, R_t \, \subset \, S$ holds. On the other hand, 
\begin{equation}\|{\rm exp}^{-1}(p)\|_{u(t)}
\ \leq \ \frac{\min \{\delta_0^a,\frac{\delta_0^b}{\alpha_0}\}}{C_1^a}
\, r\end{equation}
holds. So, in virtue of $1^\circ$ and Lemma 4.17,  

\newpage \noindent 
$$\| \, (D(\psi |_t))_{u(t)} \, ({\rm exp}^{-1}(p)) \, \|_{{\mathbb C}^n}
\ \leq \ \min \{\delta_0^a,\frac{\delta_0^b}{\alpha_0}\}
\, r$$
holds and for any $s\in[0,1]$, 
$$\| \, s \, (D(\psi |_t))_{u(t)} ({\rm exp}^{-1}(p)) \, \|_{{\mathbb C}^n}
\ < \ \delta_0^a$$
holds. Hence, in virtue of $1^\circ$ and Lemma 4.14,  
$$e^{R_t,\psi}((D(\psi |_t))_{u(t)}({\rm exp}^{-1}(p)),\psi(u(t)),t,\chi(t)) 
\ \in \ D_{\frac{r_0+1}{2}}^n,$$ 
$$\| \, e^{R_t,\psi}((D(\psi |_t))_{u(t)}({\rm exp}^{-1}(p)),\psi(u(t)),t,\chi(t)) 
\, - \, \psi(u(t)) \, \|_{{\mathbb C}^n}$$
$$= \ \| \, e^{R_t,\psi}((D(\psi |_t))_{u(t)}({\rm exp}^{-1}(p)),\psi(u(t)),t,\chi(t)) 
\, - \, e^{R_t,\psi}(0,\psi(u(t)),t,\chi(t)) \, \|_{{\mathbb C}^n}$$
$$\leq \ \alpha_0 
\, \| \, (D(\psi |_t))_{u(t)} \, ({\rm exp}^{-1}(p)) \, \|_{{\mathbb C}^n}
\ \leq \ \delta_0^b \, r$$
hold. Now, because of 
$$e^{R_t,\psi}((D(\psi |_t))_{u(t)}({\rm exp}^{-1}(p)),\psi(u(t)),t,\chi(t)) 
\ = \ \psi({\rm exp}({\rm exp}^{-1}(p))),$$
$$p\in M_\psi, \ \ \ \psi(p)\in D_{\frac{r_0+1}{2}}^n,$$
\begin{equation}\| \, \psi(p) \, - \, \psi(u(t)) \, \|_{{\mathbb C}^n} 
\ \leq \ \delta_0^b \, r\end{equation}
hold. Since from $p \, \in \, U_t \, (=U\cap\pi^{-1}(\{t\}))$ 
and $1^\circ$, $\pi(p)=t=\pi(u(t))\in X^\prime_\psi$ holds, 
$(\psi(p),\pi(p)) \, \in \, D_{\frac{r_0+1}{2}}^n \times X^\prime_\psi$ holds. 
Therefore, by $d_t(q,p)<\delta\leq\delta_2$ and Lemma 4.18, 
$q\in M_\psi$ and 
\begin{equation}\|\psi(q)-\psi(p)\|_{{\mathbb C}^n}
\ \leq \ C_2 \, d_t(q,p)\end{equation}
hold. From $(4.2)$ and $(4.3)$, 
$$\| \, \psi(p) \, - \, \psi(u(t)) \, \|_{{\mathbb C}^n} 
\ < \ \delta_0^b,$$
$$\| \, \psi(q) \, - \, \psi(u(t)) \, \|_{{\mathbb C}^n} 
\ \leq \ C_2 \, d_t(q,p) \, + \, \delta_0^b \, r$$ 
$$< \ C_2 \, \delta \, + \, \delta_0^b \, r 
\ \leq \ C_2 \, \frac{\delta_0^b}{C_2} \, (1-r) \, + \, \delta_0^b \, r
\ = \ \delta_0^b$$
hold. So, for $s\in [0,1]$, as we set 
$$c(s) \ := \ (1-s)\psi(p)+s\psi(q) 
\ = \ \psi(p)+s(\psi(q)-\psi(p)),$$ 
$$\| \, c(s) \, - \, \psi(u(t)) \, \|_{{\mathbb C}^n} 
\ < \ \delta_0^b$$
holds. Therefore, in virtue of $1^\circ$, Lemma 4.14 and $(4.3)$, 
$$(c(s),\psi(u(t)),t,\chi(t)) \ \in \ V^\prime_{R_t,\psi},$$
$$\|h^{R_t,\psi}(\psi(q),\psi(u(t)),t,\chi(t))
-h^{R_t,\psi}(\psi(p),\psi(u(t)),t,\chi(t))\|_{{\mathbb C}^n}$$
$$\leq \ \beta_0 \, \|\psi(q)-\psi(p)\|_{{\mathbb C}^n} 
\ \leq \ \beta_0 \, C_2 \, d_t(q,p)$$
hold. Hence, further, from $p \, \in \, U_t \, = \, U\cap\pi^{-1}(\{t\})$, 
$$h^{R_t,\psi}(\psi(p),\psi(u(t)),t,\chi(t))
\ = \ (D(\psi |_t))_{u(t)} ({\rm exp}^{-1}(p))$$
holds. Now, we set 
$$\dot{z} \ := \ h^{R_t,\psi}(\psi(q),\psi(u(t)),t,\chi(t)).$$
Then, 
$$\| \, \dot{z} \, - \, (D(\psi |_t))_{u(t)}({\rm exp}^{-1}(p)) \, \|_{{\mathbb C}^n}
\ \leq \ \beta_0 \, C_2 \, d_t(q,p)$$
holds. In virtue of $1^\circ$ and Lemma 4.17, 
\begin{equation}\| \, (D(\psi |_t))_{u(t)}^{-1}(\dot{z}) 
\, - \, {\rm exp}^{-1}(p) \, \|_{u(t)}
\ \leq \ \beta_0 \, C_1^b \, C_2 \, d_t(q,p)\end{equation} 
holds. From $(4.1)$ and $(4.4)$, 
$$\| \, (D(\psi |_t))_{u(t)}^{-1}(\dot{z}) \, \|_{u(t)}$$
$$\leq \ \beta_0 \, C_1^b \, C_2 \, d_t(q,p)
\, + \, \frac{\min \{\delta_0^a,\frac{\delta_0^b}{\alpha_0}\}}{C_1^a}
\, r$$
$$< \ \beta_0 \, C_1^b \, C_2 \, \delta
\, + \, \frac{\min \{\delta_0^a,\frac{\delta_0^b}{\alpha_0}\}}{C_1^a}
\, r$$
$$\leq \ \beta_0 \, C_1^b \, C_2 \, 
\frac{\min\{\delta_0^a,\frac{\delta_0^b}{\alpha_0}\}}{\beta_0C_1^aC_1^bC_2}
\, (1-r) \, + \, \frac{\min \{\delta_0^a,\frac{\delta_0^b}{\alpha_0}\}}{C_1^a}
\, r$$
$$= \ \frac{\min \{\delta_0^a,\frac{\delta_0^b}{\alpha_0}\}}{C_1^a}$$
holds. That is, 
$$(D(\psi |_t))_{u(t)}^{-1}(\dot{z}) \ \in \ B^\prime_{u(t)}(T_{\mathbb R})$$
holds. Therefore, 

\newpage \noindent 
$$q \ = \ {\rm exp}((D(\psi |_t))_{u(t)}^{-1}(\dot{z})) 
\ \in \ U\cap\pi^{-1}(\{t\}) 
\ = \ U_t$$
holds. Further, from $(4.4)$, 
$$\|\Psi(q)-\Psi(p)\|_{u(t)}$$
$$= \ \frac{C_1^a}{\min \{\delta_0^a,\frac{\delta_0^b}{\alpha_0}\}}
\, \| \, (D(\psi |_t))_{u(t)}^{-1}(\dot{z}) 
\, - \, {\rm exp}^{-1}(p) \, \|_{u(t)}$$
$$\leq \ 
\frac{\beta_0 C_1^a C_1^b C_2}{\min \{\delta_0^a,\frac{\delta_0^b}{\alpha_0}\}}
\, d_t(q,p)$$
$$< \ 
\frac{\beta_0 C_1^a C_1^b C_2}{\min \{\delta_0^a,\frac{\delta_0^b}{\alpha_0}\}}
\, \delta$$
$$\leq \ \frac{\beta_0 C_1^a C_1^b C_2}{\min \{\delta_0^a,\frac{\delta_0^b}{\alpha_0}\}}
\, \frac{ \min \{\delta_0^a,\frac{\delta_0^b}{\alpha_0}\} }{\beta_0C_1^aC_1^bC_2} \, 
\varepsilon \ = \ \varepsilon$$
holds. 
\hfill 
$\blacksquare$

\noindent 
{\sf Problem} : 

When dose there exist a global holomorphic linear connection 
(or a global holomorphic spray) on a complex manifold ? 
For example, is a closed Riemann surface 
admitting a global holomorphic linear connection an elliptic curve ? 
Also, related to this, dose there exist a nontrivial example 
of (sheaves of) manifolds 
modeled on (locally free sheaves of) the modules 
of all holomorphic sections of holomorphic vector bundles ? 
$\mbox{\ \ \ }$   \hfill ---

\noindent  
{\sf Remark} : 

Let $\Lambda$ be a topological space. 
Let $N$ be an open subset of ${\mathbb R}^m\times \Lambda$. 
Let $K$ be a compact subset of $N$. 
Then, there exists $\delta>0$ such that 
for any $x, y \, \in \, {\mathbb R}^m$ and $\lambda \in \Lambda$, 
$$(x,\lambda) \, \in \, K, \ \ \ \|y-x\|_{{\mathbb R}^m} \, < \, \delta$$
$$\Longrightarrow \ \ \ \ \ \ (y,\lambda) \, \in \, N$$
holds. 

\noindent  
{\sf Proof} : 
(We do not use this remark, but we use argument under this proof 
to prove other lemmas. 
As a preliminary announcement, we describe it here.) 

For $(x,\lambda)\in K$, there exist $\delta_{x,\lambda}>0$ 
and an open set $\Lambda_{x,\lambda}$ of $\Lambda$ such that 
$\lambda \, \in \, \Lambda_{x,\lambda}$ holds and 
for any $y \in {\mathbb R}^m$ and $\mu \in \Lambda$, 
$$\|y-x\|_{{\mathbb R}^m}<2\delta_{x,\lambda}, 
\ \mu \in \Lambda_{x,\lambda}$$
$$\Longrightarrow \ \ \ (y,\mu) \, \in \, N$$
holds. Then, there exists a finite subset $L$ of $K$ such that 
$$(x,\lambda)\in K 
\ \ \ \Longrightarrow \ \ \ 
\exists \ (x^\prime,\lambda^\prime)\in L 
\ : \ \|x-x^\prime\|_{{\mathbb R}^m}<\delta_{x^\prime,\lambda^\prime}, 
\ \lambda \in \Lambda_{x^\prime,\lambda^\prime}$$
holds. Now, we set 
$\delta \, := \, \min
(\{\delta_{x,\lambda}\}_{(x,\lambda)\in L}\cup\{1\}).$
$\delta > 0$ holds. Let 
$(x,\lambda) \, \in \, K$ and  
$\|y-x\|_{{\mathbb R}^m} \, < \, \delta$. 
Then, there exists $(x^\prime,\lambda^\prime)\in L$ such that  
$$\|x-x^\prime\|_{{\mathbb R}^m}<\delta_{x^\prime,\lambda^\prime}, 
\ \ \ \ \ \ \lambda \in \Lambda_{x^\prime,\lambda^\prime}$$ 
hold. However, 
$$\|y-x^\prime\|_{{\mathbb R}^m} 
\, \leq \, \|y-x\|_{{\mathbb R}^m}
+\|x-x^\prime\|_{{\mathbb R}^m}$$
$$< \, \delta+\delta_{x^\prime,\lambda^\prime} 
\, \leq \, 2\delta_{x^\prime,\lambda^\prime}$$
holds. $(y,\lambda) \in N$ holds. 
\hfill 
$\blacksquare$

\newpage

\sect{Proof of Lemma 1.8} 
\ \ \ \ \ $M$ is normal and $\{M_\varphi\}_{\varphi\in S}$ 
is a finite open covering of $M$. Hence, 
there exists an open covering $\{G_\varphi\}_{\varphi\in S}$ of $M$ 
such that for any $\varphi\in S$, 
$$\overline{G_\varphi} \ \subset \ M_\varphi$$
holds. For $\varphi\in S$, we set 
$$X^\prime_\varphi \ := \ \pi(G_\varphi).$$ 
Then, 
$$\overline{X^\prime_\varphi} \, = \,
\overline{\pi(G_\varphi)} \, = \, 
\pi(\overline{G_\varphi}) \ \subset \ 
\pi(M_\varphi) \, = \, X_\varphi$$
holds. In addition, since $G_\varphi$ is an open set of $M_\varphi$, 
$X^\prime_\varphi$ is an open set of $X_\varphi$. So, 
$X^\prime_\varphi$ is an open set of $X$. 
On the other hand, since $\cup_{\, \varphi \in S} \, \varphi(\overline{G_\varphi})$ 
is a compact subset of $D^n \, (:=D_1^n)$, there exists 
$r_0 \, \in \, (0,1)$ such that 
$$\cup_{\, \varphi \in S} \, \varphi(\overline{G_\varphi}) 
\ \subset \ D_{r_0}^n$$
holds. Now, we set 
$$M^\prime_\varphi \ := \ 
(\varphi,\pi |_\varphi)^{-1} \, 
(D_{r_0}^n\times X^\prime_\varphi).$$
$M^\prime_\varphi$ is an open set of $M$. 
Also, $\overline{D_{r_0}^n}\times\overline{X^\prime_\varphi}$ 
is a compact subset of $D^n\times X_\varphi$, 
$$\overline{M^\prime_\varphi} \ = \ 
(\varphi,\pi |_\varphi)^{-1} \, 
(\overline{D_{r_0}^n}\times\overline{X^\prime_\varphi})
\ \subset \ M_\varphi$$
holds. From 
$$(\varphi,\pi |_\varphi) \, (G_\varphi) 
\ \subset \ 
\varphi(G_\varphi) \, \times \, \pi(G_\varphi) 
\ \subset \ 
D_{r_0}^n \, \times \, X^\prime_\varphi,$$ 
$G_\varphi \, \subset \, M^\prime_\varphi$ holds. 
Therefore, $M \, = \, \cup_{\varphi \in S} \, G_\varphi 
\, = \, \cup_{\varphi \in S} \, M^\prime_\varphi$ holds. 
$(\, \{\, (M^\prime_\varphi,X^\prime_\varphi)\, \}_{\varphi\in S}, \, r_0 \, )$ 
is a system of compact coordinate neighborhoods of $M$. 

Similarly, there exist $\{\, (M^{\prime\prime}_\varphi,
X^{\prime\prime}_\varphi)\, \}_{\varphi\in S}$ 
and $r_1 \, \in \, (0,r_0)$ such that 
$M^{\prime\prime}_\varphi$ is an open set of $M$,  
$X^{\prime\prime}_\varphi$ is an open set of $X$ 
and 
$$\overline{M^{\prime\prime}_\varphi}\subset M^\prime_\varphi, \ \ \ 
\overline{X^{\prime\prime}_\varphi}\subset X^\prime_\varphi,$$ 
$$M^{\prime\prime}_\varphi \
= \ (\varphi, \pi |_\varphi)^{-1}
\, (D_{r_1}^n \times X^{\prime\prime}_\varphi),$$ 
$$M \, = \, \cup_{\varphi\in S} \, M^{\prime\prime}_\varphi$$
hold. Further, similarly, there exist $\{\, (M^{\prime\prime\prime}_\varphi,
X^{\prime\prime\prime}_\varphi)\, \}_{\varphi\in S}$ 
and $r_2 \, \in \, (0,r_1)$ 
such that 
$M^{\prime\prime\prime}_\varphi$ is an open set of $M$, 
$X^{\prime\prime\prime}_\varphi$ is an open set of $X$ 
and 
$$\overline{M^{\prime\prime\prime}_\varphi}\subset M^{\prime\prime}_\varphi, \ \ \ 
\overline{X^{\prime\prime\prime}_\varphi}\subset X^{\prime\prime}_\varphi,$$ 
$$M^{\prime\prime\prime}_\varphi \
= \ (\varphi, \pi |_\varphi)^{-1}
\, (D_{r_2}^n \times X^{\prime\prime\prime}_\varphi),$$ 
$$M \, = \, \cup_{\varphi\in S} \, M^{\prime\prime\prime}_\varphi$$
hold.

Now, since $X$ is normal, for $\varphi\in S$, 
there exists a non-negative continuous function $\rho_{1,\varphi}$ on $X$ 
such that $$\overline{X^{\prime\prime\prime}_\varphi} 
\ \subset \ \{ \, t\in X \, | \, \rho_{1,\varphi}(t)\not=0 \, \}
\ \subset \ X^{\prime\prime}_\varphi$$ 
holds. Also, there exists a non-negative $C^\infty$-function $\rho_2$ 
on ${\mathbb C} \, (={\mathbb C}^1={\mathbb R}^2)$ such that 
$$\{ \, (a,b)\in{\mathbb R}^2 \, | \, \rho_2(a,b)\not=0 \, \}
\ = \ D_{r_1}^1$$
holds. Now, we define a non-negative function $\rho_\varphi$ on $M$ as 
$$\rho_\varphi(p)
\ := \ 
\left\{\begin{array}{ccr}
\rho_{1,\varphi}(\pi(p)) \, (\prod_{ \, k \, = \, 1,2,\cdots,n} \, 
\rho_2(\varphi_{(k)}(p)))
& & (p\in M_\varphi),\\ 
0 & & (p\in M\setminus M_\varphi). 
\end{array}\right.$$
Here, $\varphi_{(k)}$ is the $k$-th component of $\varphi$. Now, because of 
$$\{ \, p\in M \, | \, \rho_\varphi(p)\not=0 \, \}
\ \subset \ \{ \, p\in M_\varphi 
\, | \, (\varphi,\pi |_\varphi)(p) 
\in D_{r_1}^n \times X^{\prime\prime}_\varphi \, \}
\, = \, M^{\prime\prime}_\varphi,$$
$\rho_\varphi$ is a continuous function such that 
$${\rm supp} \, (\rho_\varphi) 
\ \subset \ \overline{M^{\prime\prime}_\varphi}
\ \subset \ M^\prime_\varphi$$
holds. However, because 
$$M^{\prime\prime\prime}_\varphi
\, = \, \{ \, p\in M_\varphi 
\, | \, (\varphi,\pi |_\varphi)(p) 
\in D_{r_2}^n \times X^{\prime\prime\prime}_\varphi \, \}
\ \subset \ \{ \, p\in M \, | \, \rho_\varphi(p)\not=0 \, \}$$
also holds, 
$$M \, = \, \cup_{\varphi\in S} \, M^{\prime\prime\prime}_\varphi
\, = \, \cup_{\varphi\in S} \, \{ \, p\in M \, | \, \rho_\varphi(p)\not=0 \, \}$$
holds. 
\hfill 
$\blacksquare$

\newpage

\sect{Proof of Lemma 1.11} 
\ \ \ \ \ From the induction, it suffices to show it in the case of $|\gamma|=1$. 
Furthermore, when $|\gamma|=1$ holds, 
even if we assume $n=1$, we do not lose generality. 
Let $n=1$.

Let $(w,s)\in U$ and $\varepsilon>0$. 
Then, there exist $r>0$ and an open set $T^\prime$ of $T$ such that 
$s\in T^\prime$ and 
$$|z-w|<2r, \ t\in T^\prime \ \ \  
\Longrightarrow \ \ \ (z,t)\in U$$
hold. Further, for any $\tau\in[0,2\pi]$, 
there exist $\delta_\tau\in (0,r)$ and an open set $T_\tau$ of $T^\prime$ 
such that $s\in T_\tau$ and 
$$|\theta-\tau|<\delta_\tau, \ |z-w|<\delta_\tau, \ t\in T_\tau$$
$$\Longrightarrow$$
$$| \, f(z+re^{\sqrt{-1}\theta},t) \, 
- \, f(w+re^{\sqrt{-1}\tau},s) \, | \ 
< \ \frac{1}{2}r \varepsilon$$
hold. There exists a finite subset $K \, (\not=\emptyset)$ of $[0,2\pi]$ 
such that 
$$\theta\in[0,2\pi] \ \ \ \Longrightarrow \ \ \ 
\exists \, \tau \in K: \, |\theta-\tau|<\delta_\tau$$
holds. Now, we set 
$$\delta\, :=\, \min_{\tau\in K}\delta_\tau, \ \ \ 
T^{\prime\prime}\, :=\, \cap_{\tau\in K}T_\tau.$$
$\delta>0$ and $s\in T^{\prime\prime}$ hold. 
$T^{\prime\prime}$ is an open set of $T$. 

Now, let $|z-w|<\delta$ and $t\in T^{\prime\prime}$. Then, 
for any $\theta\in [0,2\pi]$, there exists $\tau\in K$ such that 
$$| \, f(z+re^{\sqrt{-1}\theta},t) \, - \, f(w+re^{\sqrt{-1}\theta},s) \, |$$
$$\leq \ | \, f(z+re^{\sqrt{-1}\theta},t) \, - \, f(w+re^{\sqrt{-1}\tau},s) \, |
\ + \ | \, f(w+re^{\sqrt{-1}\tau},s) \, - \, f(w+re^{\sqrt{-1}\theta},s) \, |$$
$$< \ r\varepsilon$$
holds. So, 
$$| \, (\partial_z f)(z,t) \, - \, (\partial_z f)(w,s) \, |$$
$$= \ | \, \frac{1}{2\pi r} \,  
\int_{\theta\in[0,2\pi]} \, e^{-\sqrt{-1}\theta} \, 
(\, f(z+re^{\sqrt{-1}\theta},t)
\, - \, f(w+re^{\sqrt{-1}\theta},s) \, ) 
\, d\theta \, |$$
$$< \, \varepsilon$$
holds. 
\hfill 
$\blacksquare$

\newpage

\sect{Proof of Lemma 1.16} 
\ \ \ \ \ From Lemma 1.14, 
there exist an open set $U^\prime$ of $M_\varphi$ and $C>0$ 
such that the followings hold.  

[ \ (1) \ $K\subset U^\prime$ holds. 
\ (2) \ 
For any $p\in U^\prime$ and $\dot{p}\in (T_{\mathbb R})_p(\pi^{-1}(\{\pi(p)\}))$, 
$\|\dot{p}\|_{\varphi,p}
\, \leq \, C \, \|\dot{p}\|_p$
holds. \ ] 

\noindent 
Then, $K \subset U^\prime\cap N$ holds, $U^\prime\cap N$ is an open set of $M$ 
and $K$ is a closed set of $M$. Because $M$ is normal, 
there exists an open set $U^{\prime\prime}$ of $M$ such that 
$$K \ \subset \ U^{\prime\prime} 
\ \subset \ \overline{U^{\prime\prime}}
\ \subset \ U^\prime\cap N$$
holds. Also, since $U^\prime\cap N \ \subset \ M_\varphi$ holds, 
$$K \ \subset \ U^{\prime\prime} 
\ \subset \ \overline{U^{\prime\prime}}
\ \subset \ U^\prime\cap N 
\ \subset \ M_\varphi$$
holds. So, for $p\in K$, there exist $\delta_p>0$ 
and an open set $T_p$ of $X$ such that 
$$\pi(p) \, \in \, T_p,$$
$$\{ \, z \in {\mathbb C}^n \, | 
\, \|z-\varphi(p)\|_{{\mathbb C}^n} < (C+1)\delta_p \, \}
\times T_p
\ \subset \ (\varphi,\pi |_\varphi)(U^{\prime\prime})$$
hold. Then, there exists a finite subset $L$ of $K$ such that 
$$(\varphi,\pi |_\varphi)(K)
\ \subset \ 
\cup_{p\in L}(\{ \, z \in {\mathbb C}^n \, | 
\, \|z-\varphi(p)\|_{{\mathbb C}^n} < \delta_p \, \}
\times T_p)$$
holds. Now, we set 
$$\delta:=\min(\{\delta_p\}_{p\in L}\cup\{1\}),$$
$$V:=\cup_{p\in L}(\{ \, z \in {\mathbb C}^n \, | 
\, \|z-\varphi(p)\|_{{\mathbb C}^n} < \delta_p \, \}
\times T_p),$$
$$U:=(\varphi,\pi |_\varphi)^{-1}(V).$$ 
$\delta \, \in \, (0,1]$ holds. 
$U$ is an open set of $N$. $U$ includes $K$.

Now, let $p \, \in \, U$, $q \, \in \, \pi^{-1}(\{\pi(p)\})$ 
and $d_{\pi(p)}(q,p) \, < \, \delta$. 
Then, there exists $p^\prime \in L$ such that 
$$\|\varphi(p)-\varphi(p^\prime)\|_{{\mathbb C}^n}<\delta_{p^\prime}, 
\ \ \ \pi(p) \in T_{p^\prime}$$
hold. So, for any $z\in {\mathbb C}^n$, 
$$\|z-\varphi(p)\|_{{\mathbb C}^n}<C\delta$$

\newpage \noindent  
$$\Longrightarrow$$
$$\|z-\varphi(p^\prime)\|_{{\mathbb C}^n}
\ \leq \ \|z-\varphi(p)\|_{{\mathbb C}^n}
+\|\varphi(p)-\varphi(p^\prime)\|_{{\mathbb C}^n}$$
$$<C\delta+\delta_{p^\prime}
\ \leq \ (C+1)\delta_{p^\prime}$$
$$\Longrightarrow$$
$$(z,\pi(p)) 
\ \subset \ (\varphi,\pi |_\varphi)(U^{\prime\prime})$$
holds. Therefore, 
$$ \forall \, z\in{\mathbb C}^n 
\ : \ [ \ \|z-\varphi(p)\|_{{\mathbb C}^n}<C\delta $$
$$\Longrightarrow \ [ \ z\in D^n, 
\ \varphi |_{\pi(p)}^{-1}(z)\in U^{\prime\prime} \ ] \ ]$$
holds. Now, we set 
$$H \ := \ \{ \, z \in {\mathbb C}^n 
\, | \, \|z-\varphi(p)\|_{{\mathbb C}^n}<C\delta \, \}.$$
$H$ is an open set of $D^n$. Further, we set 
$$G \ := \ \varphi |_{\pi(p)}^{-1}(H),$$ 
$$F \ := \ \pi^{-1}(\{\pi(p)\}) \, \setminus \, G.$$
$G$ is an open set of $M_{\varphi,\pi(p)}$. 
$F$ is a closed set of $\pi^{-1}(\{\pi(p)\})$. 
$$p \ \in \ G \ \subset \ U^{\prime\prime}$$
holds. Now, because of $d_{\pi(p)}(q,p) \, < \, \delta \leq 1$, 
there exists $\{c_m\}_{m \in {\mathbb N}}$ such that the followings hold. 

[ \ (1) \ For any $m \in {\mathbb N}$, $c_m$ 
is a piecewise $C^\infty$-map from $[0,1]$ to $\pi^{-1}(\{\pi(p)\})$ 
and $$c_m(0)=p, \ \ \ c_m(1)=q, \ \ \ L_{\pi(p)}(c_m)<\delta$$ 
hold. \ (2) \ $d_{\pi(p)}(q,p) 
\, = \, \lim_{m\rightarrow\infty} \, L_{\pi(p)}(c_m)$ holds. \ ] 

\noindent 
Now, for $m\in {\mathbb N}$, we set 
$$s_m \ := \ \min 
\, ( \, c_m^{-1}(F) 
\, \cup \, \{1\} \, ).$$
$s_m \, \in \, (0,1]$ holds. However, 
because of $$s\in [0,s_m) 
\ \ \ \Longrightarrow \ \ \ 
c_m(s) \in G,$$ 

\newpage \noindent 
$$c_m([0,s_m]) \ \subset \ \overline{G} 
\ \subset \ \overline{U^{\prime\prime}}
\ \subset \ U^\prime\cap N 
\ \subset \ M_\varphi$$
holds. So, because of 
$$\|\varphi(c_m(s_m))-\varphi(p)\|_{{\mathbb C}^n} 
\ \leq \ 
\int_0^{s_m} \, \| \frac{d}{ds} 
( (\varphi\circ c_m)(s) ) 
\|_{{\mathbb C}^n} \, ds
$$
$$\leq \ 
\int_0^{s_m} \, C\| \frac{d}{ds} 
( c_m(s) ) 
\|_{c_m(s)} \, ds
\ = \ CL_{\pi(p)}(\{c_m(s)\}_{s\in[0,s_m]}) 
\ < \ C\delta,$$
$\varphi(c_m(s_m)) \, \in \, H$ holds. 
Hence, $c_m(s_m) \, \in \, G$ and $s_m \, = \, 1$ hold. 
Therefore, 
$$q=c_1(1)=c_1(s_1) \ \in \ N,$$
$$\|\varphi(q)-\varphi(p)\|_{{\mathbb C}^n}$$ 
$$= \ \lim_{m\rightarrow\infty} \, \|\varphi(c_m(1))-\varphi(p)\|_{{\mathbb C}^n}$$ 
$$= \ \lim_{m\rightarrow\infty} \, \|\varphi(c_m(s_m))-\varphi(p)\|_{{\mathbb C}^n}$$ 
$$\leq \ \lim_{m\rightarrow\infty} \, CL_{\pi(p)}(\{c_m(s)\}_{s\in[0,s_m]})$$
$$= \ \lim_{m\rightarrow\infty} \, CL_{\pi(p)}(c_m) 
\ = \ C d_{\pi(p)}(q,p)$$
hold. 
\hfill 
$\blacksquare$

\newpage

\sect{Proof of Remark on Lemma 1.16} 
\ \ \ \ \ From Lemma 1.14, there exist an open set $U^\prime$ of $M_\varphi$ 
and $C>0$ such that the followings hold.  

[ \ (1) \ $K\subset U^\prime$ holds. 
\ (2) \  
For any $p\in U^\prime$ and $\dot{p}\in (T_{\mathbb R})_p(\pi^{-1}(\{\pi(p)\}))$, 
$\|\dot{p}\|_p
\, \leq \, C \, \|\dot{p}\|_{\varphi,p}$
holds. \ ] 

\noindent 
Then, $U^\prime\cap N$ is an open set of $M_\varphi$. 
$U^\prime\cap N$ includes $K$. So, for $p\in K$, 
there exist $\delta_p>0$ and an open set $T_p$ of $X$ such that 
$$\pi(p) \, \in \, T_p,$$
$$\{ \, z \in {\mathbb C}^n \, | 
\, \|z-\varphi(p)\|_{{\mathbb C}^n} < 2\delta_p \, \}
\times T_p
\ \subset \ (\varphi,\pi |_\varphi)(U^\prime\cap N)$$
hold. There exists a finite subset $L$ of $K$ such that 
$$(\varphi,\pi |_\varphi)(K)
\ \subset \ 
\cup_{p\in L}(\{ \, z \in {\mathbb C}^n \, | 
\, \|z-\varphi(p)\|_{{\mathbb C}^n} < \delta_p \, \}
\times T_p)$$
holds. Now, we set 
$$\delta:=\min(\{\delta_p\}_{p\in L}\cup\{1\}),$$
$$V:=\cup_{p\in L}(\{ \, z \in {\mathbb C}^n \, | 
\, \|z-\varphi(p)\|_{{\mathbb C}^n} < \delta_p \, \}
\times T_p),$$
$$U:=(\varphi,\pi |_\varphi)^{-1}(V).$$
$\delta>0$ holds. $U$ is an open set of $N$. $U$ includes $K$. 

Now, let $p \, \in \, U$, 
$q \, \in \, M_{\varphi,\pi(p)}$ and 
$\|\varphi(q)-\varphi(p)\|_{{\mathbb C}^n} \, < \, \delta.$ 
For $s\in[0,1]$, we set 
$$c(s):=(1-s)\varphi(p)+s\varphi(q)
=\varphi(p)+s(\varphi(q)-\varphi(p)).$$ 
Also, there exists $p^\prime \in L$ such that 
$$\|\varphi(p)-\varphi(p^\prime)\|_{{\mathbb C}^n}<\delta_{p^\prime}, 
\ \ \ \pi(p) \in T_{p^\prime}$$
hold. Hence, because of 
$$\|c(s)-\varphi(p^\prime)\|_{{\mathbb C}^n}
\leq \|\varphi(p)-\varphi(p^\prime)\|_{{\mathbb C}^n}
+s\|\varphi(q)-\varphi(p)\|_{{\mathbb C}^n}$$
$$<\delta_{p^\prime}+s\delta \leq 2\delta_{p^\prime},$$
$$(c(s),\pi(p)) \ \in \ (\varphi,\pi |_\varphi)(U^\prime\cap N)$$
holds. In particular, 
$q=(\varphi,\pi |_\varphi)^{-1}(c(1),\pi(p)) 
\, \in \, N$
holds. Further, because of 
$(\varphi,\pi |_\varphi)^{-1}(c(s),\pi(p)) 
\, \in \, U^\prime,$
$$d_{\pi(p)}(q,p)
\leq L_{\pi(p)}
(\{(\varphi,\pi |_\varphi)^{-1}(c(s),\pi(p))\}_{s\in [0,1]})$$
$$\leq \int_0^1 C\|\varphi(q)-\varphi(p)\|_{{\mathbb C}^n} ds
=C\|\varphi(q)-\varphi(p)\|_{{\mathbb C}^n}$$
holds. 
\hfill 
$\blacksquare$

\newpage

\sect{Proof of Lemma 3.7} 
\ \ \ \ \ Taking an orthonormal basis of $X$ and $Y$, respectively, 
we consider that $X=Y={\mathbb C}^n$ holds. 
$\|z\|_X=\|z\|_{{\mathbb C}^n}$ and $\|w\|_Y=\|w\|_{{\mathbb C}^n}$ holds. 
$f_k$ denotes the $k$-th component of $f$. 
$\{e_k\}_{k=1}^n$ denotes the canonical ${\mathbb Z}$-basis of ${\mathbb Z}^n$, 
that is, $\{e_k\}_{k=1}^n$ is the identity matrix of the order $n$. 

Now, as $z\in {\mathbb C}^n$ and  
$\|z\|_{{\mathbb C}^n}\leq \frac{1}{2}\varepsilon$ hold, 
$|\frac{\partial f_i}{\partial z_j}(z)|
\, = \, |\frac{2}{\pi \varepsilon}
\int_{\theta\in[0,2\pi]}e^{-\sqrt{-1}\theta} 
f_i(z+\frac{1}{4}\varepsilon e^{\sqrt{-1}\theta}e_j) d\theta| 
\, \leq \, \frac{4}{\varepsilon}$
holds. So, 
$$\|z\|_{{\mathbb C}^n}\leq \frac{1}{2}\varepsilon 
\ \ \ \Longrightarrow \ \ \ 
|\frac{\partial f_i}{\partial z_j}(z)| \ \leq \ \frac{4}{\varepsilon}$$
holds. Similarly, as $z\in {\mathbb C}^n$, 
$\|z\|_{{\mathbb C}^n}\leq \frac{1}{4}\varepsilon$  
and $s\in [0,1]$ hold, 
$|\frac{\partial^2 f_i}{\partial z_j\partial z_k}(sz)|
\, = \, |\frac{2}{\pi \varepsilon}
\int_{\theta\in[0,2\pi]}e^{-\sqrt{-1}\theta} 
\frac{\partial f_i}{\partial z_k}
(sz+\frac{1}{4}\varepsilon e^{\sqrt{-1}\theta}e_j) d\theta| 
\, \leq \, \frac{16}{\varepsilon^2}$ 
and 
$|f_i(z)-(f_i(0)+(Df_i)_0(z))|
\, = \, |\sum_{j,k} (\int_0^1 
( \frac{\partial^2 f_i}{\partial z_j\partial z_k} (sz) )
(1-s) 
ds)z_jz_k|
\, \leq \, \frac{16n^2}{\varepsilon^2}\|z\|_{{\mathbb C}^n}^2$
hold. So, 
$$\|z\|_{{\mathbb C}^n}\leq \frac{1}{4}\varepsilon 
\ \ \ \Longrightarrow \ \ \ 
|f_i(z)-(f_i(0)+(Df_i)_0(z))|\leq
\frac{16n^2}{\varepsilon^2}\|z\|_{{\mathbb C}^n}^2$$
holds. 
\hfill 
$\blacksquare$

\newpage

\sect{Proof of Lemma 4.9} 
\ \ \ \ \ Proofs of this type of lemma are often done 
by showing that the solutions can be extended 
using Gronwall inequalities, 
but here we see that the principle of contraction mappings 
works directly on an appropriate integral equation. 

For any $s\in [0,1]$, $(u_0(s),\lambda_0)\in U$ holds. The map 
$$s \, \in \, [0,1] \ \mapsto \ 
(D_xf)(u_0(s),\lambda_0) \, \in  \, L({\mathbb R}^m;{\mathbb R}^m)$$
is continuous. Hence, there exists 
$A \, : \, [0,1] \, \rightarrow \, L({\mathbb R}^m;{\mathbb R}^m)$ 
such that 
$$\frac{d}{ds}A(s) \ = \ ((D_xf)(u_0(s),\lambda_0)) \, A(s), 
\ \ \ \ \ \ A(0) \ = \ 1_{{\mathbb R}^m}$$
hold. Further, $\det(A(s))\not=0$ and 
$$1 \ \leq \ C_1 \ := 
\ \sup_{s\in[0,1]}\|(A(s))^{-1}\|_{L({\mathbb R}^m;{\mathbb R}^m)} 
\, < \, +\infty,$$ 
$$1 \ \leq \ C_2 \ := 
\ \sup_{s\in[0,1]}\|A(s)\|_{L({\mathbb R}^m;{\mathbb R}^m)} 
\, < \, +\infty$$
hold. Then, for $s\in [0,1]$, there exist $\delta_s>0$ 
and an open set $\Lambda_s$ of $\Lambda$ such that 
$\lambda_0\in\Lambda_s$ and 
$$s^\prime \, \in \, [0,1], 
\ |s^\prime -s| \, < \, \delta_s, 
\ x \, \in \, {\mathbb R}^m, 
\ \|x\|_{{\mathbb R}^m} \, < \, 2\delta_s, 
\ \lambda \, \in \, \Lambda_s$$
$$\Longrightarrow$$
$$\left\{\begin{array}{c} 
(u_0(s^\prime)+A(s^\prime)x,\lambda) \, \in \, U, \\ 
\|(D_xf)(u_0(s^\prime)+A(s^\prime)x,\lambda)
-(D_xf)(u_0(s),\lambda_0)\|_{L({\mathbb R}^m;{\mathbb R}^m)} 
\ \leq \ \frac{1}{4C_1C_2}\end{array}\right.$$
hold. There exists a finite subset $L \, (\not=\emptyset)$ of $[0,1]$ such that 
$$s\in[0,1] \ \ \ \Longrightarrow 
\ \ \ \exists \, s^\prime \, \in \, L \, : \, 
|s-s^\prime|<\delta_{s^\prime}$$
holds. Now, we set 
$$\delta^\prime \ := \ \min_{s\in L} \, \delta_s,$$
$$\Lambda^\prime \ := \ \cap_{s\in L} \, \Lambda_s.$$
$\delta^\prime>0$ and $\lambda_0\in \Lambda^\prime$ hold. 
$\Lambda^\prime$ is an open set of $\Lambda$. 
Now, as $s\in [0,1]$, 
$x\in{\mathbb R}^m$, $\|x\|_{{\mathbb R}^m}\leq\delta^\prime$ 
and $\lambda \in \Lambda^\prime$ hold, because there exists $s^\prime\in L$ 
such that 

\newpage \noindent  
$$|s-s^\prime|<\delta_{s^\prime}$$ 
holds, from $\lambda\in\Lambda^\prime\subset\Lambda_{s^\prime}$ 
and $\|x\|_{{\mathbb R}^m}\leq\delta^\prime\leq\delta_{s^\prime}
<2\delta_{s^\prime}$, 
$$C_1 \, C_2 \, \|(D_xf)(u_0(s)+A(s)x,\lambda)
-(D_xf)(u_0(s),\lambda_0)\|_{L({\mathbb R}^m;{\mathbb R}^m)}$$ 
$$\leq \ C_1 \, C_2 \, \|(D_xf)(u_0(s)+A(s)x,\lambda)
-(D_xf)(u_0(s^\prime),\lambda_0)\|_{L({\mathbb R}^m;{\mathbb R}^m)}$$ 
$$+ \ C_1 \, C_2 \, \|(D_xf)(u_0(s),\lambda_0)
-(D_xf)(u_0(s^\prime),\lambda_0)\|_{L({\mathbb R}^m;{\mathbb R}^m)}$$ 
$$\ \leq \ \frac{1}{4}+\frac{1}{4} \ \leq \ \frac{1}{2}$$
holds. That is, 
\begin{equation}s\in [0,1], 
\ x\in{\mathbb R}^m, \ \|x\|_{{\mathbb R}^m}\leq\delta^\prime, 
\ \lambda \in \Lambda^\prime\end{equation}
$$\Longrightarrow$$
$$\left\{\begin{array}{c} 
(u_0(s)+A(s)x,\lambda) \, \in \, U, \\ 
C_1 \, C_2 \, \|(D_xf)(u_0(s)+A(s)x,\lambda)
-(D_xf)(u_0(s),\lambda_0)\|_{L({\mathbb R}^m;{\mathbb R}^m)}
\ \leq \ \frac{1}{2}\end{array}\right.$$ 
holds. 

Now, let $\varepsilon>0$. Then, for $s\in [0,1]$, 
there exist $\delta^{\prime\prime}_s>0$ 
and an open set $\Lambda^{\prime\prime}_s$ of $\Lambda$ such that 
$\lambda_0\in\Lambda^{\prime\prime}_s$ and 
$$s^{\prime\prime} \, \in \, [0,1], 
\ |s^{\prime\prime}-s| \, < \, \delta^{\prime\prime}_s, 
\ \lambda \, \in \, \Lambda^{\prime\prime}_s$$
$$\Longrightarrow$$
$$\left\{\begin{array}{c} 
(u_0(s^{\prime\prime}),\lambda) \, \in \, U, \\ 
\|f(u_0(s^{\prime\prime}),\lambda)
-f(u_0(s),\lambda_0)\|_{{\mathbb R}^m} 
\ \leq \ \frac{1}{8C_1}\min\{\delta^\prime,\frac{\varepsilon}{2C_2}\}
\end{array}\right.$$
hold. There exists a finite subset 
$L^{\prime\prime} \, (\not=\emptyset)$ of $[0,1]$ such that 
$$s\in[0,1] \ \ \ \Longrightarrow 
\ \ \ \exists \, s^{\prime\prime} \, \in \, L^{\prime\prime} \, : \, 
|s-s^{\prime\prime}|<\delta_{s^{\prime\prime}}^{\prime\prime}$$
holds. Now, we set 
$$\Lambda^{\prime\prime\prime} \ := \ 
\cap_{s\in L^{\prime\prime}} \, \Lambda^{\prime\prime}_s.$$
$\lambda_0\in \Lambda^{\prime\prime\prime}$ holds. 
$\Lambda^{\prime\prime\prime}$ is an open set of $\Lambda$. 
Now, as $s\in [0,1]$ and 
$\lambda \in \Lambda^{\prime\prime\prime}$ hold, because 
there exists $s^{\prime\prime}\in L^{\prime\prime}$ such that 
$$|s-s^{\prime\prime}|<\delta^{\prime\prime}_{s^{\prime\prime}}$$ 
holds, from $\lambda\in\Lambda^{\prime\prime\prime}
\subset\Lambda^{\prime\prime}_{s^{\prime\prime}}$, 

\newpage \noindent  
$$C_1 \, \|f(u_0(s),\lambda)-f(u_0(s),\lambda_0)\|_{{\mathbb R}^m}$$ 
$$\leq \ C_1 \, \|f(u_0(s),\lambda)
-f(u_0(s^{\prime\prime}),\lambda_0)\|_{{\mathbb R}^m}$$ 
$$+ \ C_1 \, \|f(u_0(s),\lambda_0)
-f(u_0(s^{\prime\prime}),\lambda_0)\|_{{\mathbb R}^m}$$ 
$$\ \leq \ \frac{1}{8}\min\{\delta^\prime,\frac{\varepsilon}{2C_2}\}
+\frac{1}{8}\min\{\delta^\prime,\frac{\varepsilon}{2C_2}\}
\ = \ \frac{1}{4}\min\{\delta^\prime,\frac{\varepsilon}{2C_2}\}$$
holds. That is, 
\begin{equation}s\in [0,1], 
\ \lambda \in \Lambda^{\prime\prime\prime}\end{equation}
$$\Longrightarrow$$
$$\left\{\begin{array}{c} 
(u_0(s),\lambda) \, \in \, U,\\ 
C_1 \, \|f(u_0(s),\lambda)
-f(u_0(s),\lambda_0)\|_{{\mathbb R}^m}
\ \leq \ \frac{1}{4}\min\{\delta^\prime,\frac{\varepsilon}{2C_2}\}
\end{array}\right.$$ 
holds. Now, we set 
$$V \ := \ \{ \, x\in{\mathbb R}^m 
\, | \, \|x-x_0\|_{{\mathbb R}^m} < 
\frac{1}{4} \min\{\delta^\prime,\frac{\varepsilon}{2C_2}\} \, \} 
\, \times \, (\Lambda^\prime\cap\Lambda^{\prime\prime\prime}).$$
Then, $(x_0,\lambda_0) \, \in \, V$ holds. 
Also, $V$ is an open set of ${\mathbb R}^m \, \times \, \Lambda$. 
However, from $u_0(0)=x_0$, $A(0)=1_{{\mathbb R}^m}$ and $(10.1)$, 
$V \, \subset \, U$ holds. $V$ is an open set of $U$. 

Well, we show that (2) is established. 
For $v \, \in \, C([0,1];{\mathbb R}^m)$, we set 
$$\|v\|_{C([0,1];{\mathbb R}^m)} 
\ := \ \sup_{s\in [0,1]}\|v(s)\|_{{\mathbb R}^m}.$$
We set $$\tilde{X} \ := \ 
\{ \, v \in C([0,1];{\mathbb R}^m) 
\, | \, \|v\|_{C([0,1];{\mathbb R}^m)}
\leq \min\{\delta^\prime,\frac{\varepsilon}{2C_2}\} \, \}.$$ 
$\tilde{X}$ is a complete metric space. $\tilde{X}\not=\emptyset$ holds. 
Let $(x,\lambda) \in V$. Then, from $(10.1)$, 
for any $v \, \in \, \tilde{X}$ and $s \, \in \, [0,1]$, 
$(u_0(s)+A(s)v(s),\lambda) \, \in \, U$ holds. Now, 
for $v \, \in \, \tilde{X}$ and $s \, \in \, [0,1]$, we set 
$$(F(v))(s)$$
$$:= \ (x-x_0) \, + \, \int_0^s \, (A(\sigma))^{-1}$$ 
$$( \, f(u_0(\sigma)+A(\sigma)v(\sigma),\lambda) 
\, - \, ( f(u_0(\sigma),\lambda_0) 
+ ((D_xf)(u_0(\sigma),\lambda_0)) A(\sigma) v(\sigma) ) \, ) \, d\sigma.$$
We show that $F$ is a contraction mapping of $\tilde{X}$. 
As $v\in\tilde{X}$, $s\in[0,1]$, $\sigma\in[0,1]$ and $\tau\in[0,1]$ hold, 
because of $$\| \, \tau \, v(\sigma) \, \|_{{\mathbb R}^m} 
\ \leq \ \tau \, \min\{\delta^\prime, \frac{\varepsilon}{2C_2}\} \ \leq \ \delta^\prime,$$
from $(10.1)$ and $(10.2)$, 
$$(u_0(\sigma)+\tau A(\sigma)v(\sigma),\lambda) \ \in \ U,$$
$$f(u_0(\sigma)+A(\sigma)v(\sigma),\lambda)-f(u_0(\sigma),\lambda)$$
$$= \ ( \, \int_0^1 \, (D_xf)(u_0(\sigma)+\tau A(\sigma)v(\sigma),\lambda) 
\, d\tau \, ) \, A(\sigma) \, v(\sigma),$$
$$\| \, (A(\sigma))^{-1} \, 
( \, (f(u_0(\sigma)+A(\sigma)v(\sigma),\lambda)-f(u_0(\sigma),\lambda))
\, - \, ((D_xf)(u_0(\sigma),\lambda_0)) A(\sigma) v(\sigma) \, ) \, \|_{{\mathbb R}^n}$$
$$= \ \| \, (A(\sigma))^{-1} 
\, ( \, \int_0^1 \, ((D_xf)(u_0(\sigma)+\tau A(\sigma)v(\sigma),\lambda)
-(D_xf)(u_0(\sigma),\lambda_0)) \, d\tau \, ) 
\, A(\sigma) \, v(\sigma) \, \|_{{\mathbb R}^m}$$
$$\leq \ \frac{1}{2} \, \|v(\sigma)\|_{{\mathbb R}^m}
\ \leq \ \frac{1}{2} \, \|v\|_{C([0,1];{\mathbb R}^m)},$$ 
$$\|(F(v))(s)\|_{{\mathbb R}^m}$$
$$\leq \ \|x-x_0\|_{{\mathbb R}^m} 
\ + \ \| \, \int_0^s \, (A(\sigma))^{-1} \, 
( f(u_0(\sigma),\lambda) - f(u_0(\sigma),\lambda_0) ) \, d\sigma \, \|_{{\mathbb R}^m}
\ + \ \frac{1}{2} \|v\|_{C([0,1];{\mathbb R}^m)}$$
$$\leq \ \frac{1}{4} \min\{\delta^\prime,\frac{\varepsilon}{2C_2}\}
+\frac{1}{4} \min\{\delta^\prime,\frac{\varepsilon}{2C_2}\}
+\frac{1}{2} \min\{\delta^\prime,\frac{\varepsilon}{2C_2}\}
\ = \ \min\{\delta^\prime,\frac{\varepsilon}{2C_2}\}$$
hold. So, $F(\tilde{X}) \, \subset \, \tilde{X}$ holds. 
Next, as $v\in\tilde{X}$, $v^\prime\in\tilde{X}$, 
$s\in[0,1]$, $\sigma\in[0,1]$ and $\tau\in[0,1]$ hold and we set 
$$c_{\tau} \ := \ (1-\tau)v+\tau v^\prime
\ = \ v+\tau(v^\prime-v),$$
because of 
$$\| \, c_\tau(\sigma) \, \|_{{\mathbb R}^m} 
\ \leq \ (1-\tau) \, \min\{\delta^\prime, \frac{\varepsilon}{2C_2}\} 
\, + \, \tau \, \min\{\delta^\prime, \frac{\varepsilon}{2C_2}\} 
\ \leq \ \delta^\prime,$$
from $(10.1)$, 

\newpage \noindent 
$$(u_0(\sigma)+A(\sigma)c_\tau(\sigma),\lambda) \ \in \ U,$$
$$f(u_0(\sigma)+A(\sigma)v^\prime(\sigma),\lambda)
-f(u_0(\sigma)+A(\sigma)v(\sigma),\lambda)$$
$$= \ ( \, \int_0^1 \, (D_xf)(u_0(\sigma)+A(\sigma)c_\tau(\sigma),\lambda) 
\, d\tau \, ) \, A(\sigma) \, (v^\prime(\sigma)-v(\sigma)),$$
$$\|(F(v^\prime))(s)-(F(v))(s)\|_{{\mathbb R}^m}$$
$$= \ \| \, \int_0^s \, (A(\sigma))^{-1}$$
$$( \, ( f(u_0(\sigma)+A(\sigma)v^\prime(\sigma),\lambda)
- f(u_0(\sigma)+A(\sigma)v(\sigma),\lambda) )$$
$$- \, ((D_xf)(u_0(\sigma),\lambda_0)) A(\sigma) (v^\prime(\sigma)-v(\sigma)) \, ) 
\, d\sigma \, \|_{{\mathbb R}^n}$$
$$= \ \| \, \int_0^s \, (A(\sigma))^{-1} 
\, ( \, \int_0^1 \, ((D_xf)(u_0(\sigma)+A(\sigma)c_\tau(\sigma),\lambda)
-(D_xf)(u_0(\sigma),\lambda_0)) \, d\tau \, )$$
$$A(\sigma) \, (v^\prime(\sigma)-v(\sigma)) \, d\sigma \, \|_{{\mathbb R}^m}$$
$$\leq \ \frac{1}{2} \, \|v^\prime-v\|_{C([0,1];{\mathbb R}^m)}$$ 
hold. $F$ is a contraction mapping of $\tilde{X}$. 
There exists $v\in\tilde{X}$ such that 
$$v=F(v)$$
holds. We set 
$$u \ := \ u_0+Av.$$
Then, 
$$\frac{d}{ds}u \ = \ \frac{d}{ds}u_0+(\frac{d}{ds}A)v+A(\frac{d}{ds}v)$$
$$= \ f(u_0,\lambda_0)+((D_xf)(u_0,\lambda_0))Av+A(\frac{d}{ds}F(v))$$
$$= \ f(u_0,\lambda_0)+((D_xf)(u_0,\lambda_0))Av+A A^{-1} 
( f(u,\lambda) - ( f(u_0,\lambda_0)+((D_xf)(u_0,\lambda_0))Av ) )$$
$$= \ f(u,\lambda),$$
$$u(0) \ = \ u_0(0)+A(0)v(0)
\ = \ x_0+(F(v))(0) \ = \ x_0+(x-x_0) \ = \ x,$$
$$\|u-u_0\|_{C([0,1];{\mathbb R}^m)}
\ = \ \|Av\|_{C([0,1];{\mathbb R}^m)}
\ \leq \ C_2 \, \min\{\delta^\prime,\frac{\varepsilon}{2C_2}\}
\ \leq \ \frac{\varepsilon}{2}
\ < \ \varepsilon$$
hold.  
\hfill 
$\blacksquare$

\newpage

\sect{Proof of Lemma 4.10} 
\ \ \ \ \ Let $(x_0,y_0) 
\, \in \, V \, (\subset {\mathbb C}^m={\mathbb R}^{2m}).$ 
We denote the Frechet derivative of $w$ for the variable $(x, y) \, \in \, V$ 
at a point $(s, x_0, y_0)$ by $A(s)$. 
That is, we set 
$$A(s):=(D_{(x,y)}w)(s,x_0,y_0).$$
Then, $A \, : \, [0,1] \rightarrow L({\mathbb R}^{2m};{\mathbb R}^{2m})$ 
satisfies 
$$\frac{d}{d s} A 
\ = \ ( (D_{(x,y)}f) (w(s,x_0,y_0)) ) \, A,$$
$$A(0)=1_{{\mathbb R}^{2m}}.$$
On the other hand, since $f$ is holomorphic, 
$(D_{(x,y)}f)(w(s,x_0,y_0))$ is complex linear. 
Therefore, for any $c \, \in \, {\mathbb C}$ and 
$(p,q) \, \in \, {\mathbb C}^m={\mathbb R}^{2m}$, 
both $c \, ( ( A(s) ) (p,q) )$ 
and $( A(s) ) \, (c(p,q))$ 
are the solutions of the initial value problem 
$$\frac{d}{d s}(u,v) 
\ = \ ((D_{(x,y)}f)(w(s,x_0,y_0))) \, (u,v),$$
$$(u,v)(0)=c(p,q)$$
to coincide. $(D_{(x,y)}w)(s,x_0,y_0)$ 
is complex linear. 
\hfill 
$\blacksquare$

\newpage

\sect{Proof of Lemma 4.12} 
\ \ \ \ \ This is the same as proofs of usual inverse function theorems. 
Just to be sure, we confirm it. 

When $\Theta=\emptyset$ holds, it is trivial. 
Let $\Theta \, \not= \, \emptyset$ and $\varepsilon \, > \, 0$. We set 
$$\varepsilon^\prime \ := \ \frac{\varepsilon}{1+\varepsilon}.$$
$0<\varepsilon^\prime<\min\{\varepsilon,1\}$ holds. 
For $\lambda\in\Theta$, there exist $\delta^{\prime\prime}_\lambda>0$ 
and an open set $\Lambda^{\prime\prime}_\lambda$ of $\Lambda$ such that 
$\lambda\in\Lambda^{\prime\prime}_\lambda$ and 
$$x\in{\mathbb R}^m, 
\ \|x\|_{{\mathbb R}^m}\leq\delta^{\prime\prime}_\lambda, 
\ \lambda^{\prime\prime}\in\Lambda^{\prime\prime}_\lambda$$
$$\Longrightarrow \ \ \ 
(x,\lambda^{\prime\prime})\in U, 
\ \|(D_xf)(x,\lambda^{\prime\prime})-1_{{\mathbb R}^m}
\|_{L({\mathbb R}^m;{\mathbb R}^m)}
\leq\varepsilon^\prime$$
hold. Then, there exists a finite subset 
$L \, (\not=\emptyset)$ of $\Theta$ such that 
$$\Theta \ \subset \ \Lambda^\prime 
\ := \ \cup_{\lambda\in L}\Lambda^{\prime\prime}_\lambda$$
holds. $\Lambda^\prime$ is an open set of $\Lambda$. Now, we set 
$$\delta^\prime \ := \ \min_{\lambda\in L} 
\, \delta^{\prime\prime}_\lambda,$$
$$E \ := \ \{ \, x\in{\mathbb R}^m 
\, | \, \|x\|_{{\mathbb R}^m}\leq\delta^\prime \, \}
\, \times \, \Lambda^\prime.$$
$\delta^\prime>0$ holds. Also, as $(x,\lambda)\in E$ holds, 
because there exists $\lambda^{\prime\prime}\in L$ such that 
$\lambda\in\Lambda^{\prime\prime}_{\lambda^{\prime\prime}}$ holds, 
from $ \|x\|_{{\mathbb R}^m}\leq\delta^\prime\leq
\delta^{\prime\prime}_{\lambda^{\prime\prime}}$, 
$(x,\lambda)\in U$ and  
$\|(D_xf)(x,\lambda)-1_{{\mathbb R}^m}
\|_{L({\mathbb R}^m;{\mathbb R}^m)}
\leq\varepsilon^\prime$
hold. So, 
$$\{0\} \times \Lambda^\prime \ \subset \ E \ \subset \ U,$$
\begin{equation}\sup_{(x,\lambda)\in E}\|1_{{\mathbb R}^m}-(D_xf)(x,\lambda)
\|_{L({\mathbb R}^m;{\mathbb R}^m)}
\ \leq \ \varepsilon^\prime\end{equation}
hold. Also, we set 
$$E^\prime \ :=  \ \cup_{\lambda\in\Lambda^\prime} 
\, ( \, \{ \, y\in{\mathbb R}^m 
\, | \, \|y-f(0,\lambda)\|_{{\mathbb R}^m}\leq(1-\varepsilon^\prime)\delta^\prime \, \}
\times\{\lambda\} \, ),$$
$$E^{\prime\prime} \ := \ 
\{ \, x\in{\mathbb R}^m 
\, | \, \|x\|_{{\mathbb R}^m}\leq\delta^\prime \, \} 
\, \times \, E^\prime.$$
Then, 

\newpage \noindent  
$$(x,\lambda) \in E \ \ \ \Longrightarrow 
\ \ \ (x,f(0,\lambda),\lambda) \in E^{\prime\prime},$$
$$(x,y,\lambda) \in E^{\prime\prime} 
\ \ \ \Longrightarrow \ \ \ (x,\lambda) \in E$$
hold. We define a map 
$h \, : \, E^{\prime\prime} \, \rightarrow \, {\mathbb R}^m$ as 
$$h(x,y,\lambda) \ := \ x-f(x,\lambda)+y.$$
As $(x,y,\lambda)\in E^{\prime\prime}$ and $s\in [0,1]$ hold, 
$(sx,\lambda)\in E$ and 
$$h(x,y,\lambda) \ = \ ((x-f(x,\lambda))-(0-f(0,\lambda))) \, + \, (y-f(0,\lambda))$$
$$= \ ( \, \int_0^1 \, (1_{{\mathbb R}^m}-(D_xf)(sx,\lambda)) \, ds \, ) \, x
\, + \, (y-f(0,\lambda))$$
hold. Hence, by $(12.1)$, 
\begin{equation}
(x,y,\lambda)\in E^{\prime\prime}
\end{equation}
$$\Longrightarrow$$
$$\|h(x,y,\lambda)\|_{{\mathbb R}^m} 
\ \leq \ \varepsilon^\prime \, \|x\|_{{\mathbb R}^m}
\, + \, \|y-f(0,\lambda)\|_{{\mathbb R}^m}
\ \leq \ \delta^\prime$$
holds. As $(x,y,\lambda)\in E^{\prime\prime}$, 
$(x^\prime,y,\lambda)\in E^{\prime\prime}$ and $s\in [0,1]$ hold 
and we set $c(s) \, := \, (1-s)x + sx^\prime \, = \, x+s(x^\prime-x)$, 
$(c(s),\lambda)\in E$ and 
$$h(x^\prime,y,\lambda)-h(x,y,\lambda) 
\ = \ (x^\prime-f(x^\prime,\lambda))-(x-f(x,\lambda))$$
$$= \ ( \, \int_0^1 \, (1_{{\mathbb R}^m}-(D_xf)(c(s),\lambda)) \, ds \, ) 
\, (x^\prime-x)$$
hold. Hence, by $(12.1)$, 
\begin{equation}
(x,y,\lambda)\in E^{\prime\prime}, 
\ (x^\prime,y,\lambda)\in E^{\prime\prime}
\end{equation}
$$\Longrightarrow$$
$$\|h(x^\prime,y,\lambda)-h(x,y,\lambda)\|_{{\mathbb R}^m} 
\ \leq \ \varepsilon^\prime \, \|x^\prime-x\|_{{\mathbb R}^m}$$
holds. 

For $\ell \, \in \, C(E^\prime;{\mathbb R}^m)$, we set 
$$\|\ell\|_{C(E^\prime;{\mathbb R}^m)} 
\ := \ \sup_{(y,\lambda)\in E^\prime}\|\ell(y,\lambda)\|_{{\mathbb R}^m}.$$ 
We set 
$$\tilde{X} \ := \ 
\{ \, \ell \in C(E^\prime;{\mathbb R}^m) 
\, | \, \|\ell\|_{C(E^\prime;{\mathbb R}^m)}
\leq \delta^\prime \, \}.$$ 
$\tilde{X}$ is a complete metric space. $\tilde{X}\not=\emptyset$ holds. 
For any $\ell \, \in \, \tilde{X}$ and $(y,\lambda) \, \in \, E^\prime$, 
$(\ell(y,\lambda),y,\lambda) \, \in \, E^{\prime\prime}$ holds. 
Now, for $\ell \, \in \, \tilde{X}$ and $(y,\lambda) \, \in \, E^\prime$, 
we set $$(F(\ell))(y,\lambda) 
\ := \ h(\ell(y,\lambda),y,\lambda).$$ 
As $\ell \, \in \, \tilde{X}$ holds, 
because of $F(\ell) \, \in \, C(E^\prime;{\mathbb R}^m)$, 
from $(12.2)$, $F(\ell) \, \in \, \tilde{X}$ holds. 
Further, as $\ell \, \in \, \tilde{X}$ and  
$\ell^\prime \, \in \, \tilde{X}$ hold, from $(12.3)$, 
$\|F(\ell^\prime)-F(\ell)\|_{C(E^\prime;{\mathbb R}^m)}
\ \leq \ \varepsilon^\prime \, \|\ell^\prime-\ell\|_{C(E^\prime;{\mathbb R}^m)}$
holds. $F$ is a contraction mapping of $\tilde{X}$. 
Therefore, there exists $\ell \, \in \, \tilde{X}$ such that 
$$\ell \, = \, F(\ell)$$
holds. In particular, 
$$(y,\lambda) \, \in \, E^\prime \ \ \ \ \ \ 
\Longrightarrow \ \ \ \ \ \ 
(\ell(y,\lambda),\lambda) \, \in \, E, 
\ f(\ell(y,\lambda),\lambda) \, = \, y$$
holds. 

Now, we set 
$$V^\prime \ := \ \cup_{\lambda\in\Lambda^\prime} 
\, ( \, \{ \, y\in{\mathbb R}^m 
\, | \, \|y-f(0,\lambda)\|_{{\mathbb R}^m} 
< (1-\varepsilon^\prime)\delta^\prime \, \}
\times\{\lambda\} \, ),$$
$$g \ := \ \ell_{\upharpoonright V^\prime},$$
$$U^\prime \ := \ \{ \, (g(y,\lambda),\lambda) 
\in E \, | \, (y,\lambda)\in V^\prime \, \}.$$ 
$V^\prime$ is an open set of ${\mathbb R}^m\times\Lambda$. 
$g \, \in \, C(V^\prime;{\mathbb R}^m)$ and 
$$(y,\lambda) \, \in \, V^\prime \ \ \ \ \ \ 
\Longrightarrow \ \ \ \ \ \ 
(g(y,\lambda),\lambda) \, \in \, U^\prime, 
\ f(g(y,\lambda),\lambda) \, = \, y$$
hold. On the other hand, as $(x,\lambda)\in U^\prime$ holds, 
because there exists $y\in {\mathbb R}^m$ such that 
$(y,\lambda)\in V^\prime$ and 
$x=g(y,\lambda)$ hold, 
$$y=f(g(y,\lambda),\lambda)=f(x,\lambda),$$
$$(f(x,\lambda),\lambda)=(y,\lambda)\in V^\prime,$$
$$x=g(y,\lambda)=g(f(x,\lambda),\lambda)$$
hold. That is, 

\newpage \noindent  
$$(x,\lambda) \, \in \, U^\prime \ \ \ \ \ \ 
\Longrightarrow \ \ \ \ \ \ 
(f(x,\lambda),\lambda) \, \in \, V^\prime, 
\ g(f(x,\lambda),\lambda) \, = \, x$$
holds. Therefore, homeomorphically 
the map $(x,\lambda) \, \mapsto \, (f(x,\lambda),\lambda)$ 
maps $U^\prime$ to $V^\prime$.

Also, from $(12.2)$, for any $(y,\lambda) \, \in \, V^\prime$, 
$\|g(y,\lambda)\|_{{\mathbb R}^m} \, = \, \|\ell(y,\lambda)\|_{{\mathbb R}^m} 
\, = \, \|(F(\ell))(y,\lambda)\|_{{\mathbb R}^m} 
\, = \, \|h(\ell(y,\lambda),y,\lambda)\|_{{\mathbb R}^m} 
\, \leq \, \varepsilon^\prime\delta^\prime+\|y-f(0,\lambda)\|_{{\mathbb R}^m} 
\, < \, \varepsilon^\prime\delta^\prime+(1-\varepsilon^\prime)\delta^\prime 
\, = \, \delta^\prime$
holds. Therefore, 
$$U^\prime \ \subset \ G \ := \ \{ \, x\in{\mathbb R}^m 
\, | \, \|x\|_{{\mathbb R}^m}<\delta^\prime \, \}
\, \times \, \Lambda^\prime
\ \subset \ E \ \subset \ U$$
holds. So, 
$$(x,\lambda) \, \in \, U^\prime \ \ \ \ \ \ 
\Longrightarrow \ \ \ \ \ \ 
(x,\lambda) \, \in \, G, 
\ (f(x,\lambda),\lambda) \, \in \, V^\prime$$
holds. Conversely, if $(x,\lambda) \, \in \, G$ 
and $(f(x,\lambda),\lambda) \, \in \, V^\prime$ hold, then 
$$g(f(x,\lambda),\lambda)
\ = \ \ell(f(x,\lambda),\lambda) 
\ = \ (F(\ell))(f(x,\lambda),\lambda)$$
$$= \ h(  \ell(f(x,\lambda),\lambda),  f(x,\lambda),  \lambda  )
\ = \ h(     g(f(x,\lambda),\lambda),  f(x,\lambda),  \lambda  ),$$
$$x \ = \  h(   x,   f(x,\lambda),   \lambda  )$$
hold and so, from $(12.3)$, 
$$\|g(f(x,\lambda),\lambda)-x\|_{{\mathbb R}^m}$$
$$\leq \ \|h(     g(f(x,\lambda),\lambda),  f(x,\lambda),  \lambda  )
-h(   x,   f(x,\lambda),   \lambda  )\|_{{\mathbb R}^m}$$
$$\leq \ \varepsilon^\prime \, \|g(f(x,\lambda),\lambda)-x\|_{{\mathbb R}^m},$$
$$g(f(x,\lambda),\lambda) \ = \ x,$$
$$(x,\lambda) \ = \ (g(f(x,\lambda),\lambda),\lambda) 
\ \in \ U^\prime$$
hold. Thus, 
\begin{equation}(x,\lambda) \, \in \, G, 
\ \ \ (f(x,\lambda),\lambda) \, \in \, V^\prime\end{equation}
$$\Longleftrightarrow \ \ \ \ \ \ (x,\lambda) \, \in \, U^\prime$$
holds. That is, 
$U^\prime \, = \, \{ \, (x,\lambda)\in G 
\, | \, (f(x,\lambda),\lambda)\in V^\prime \, \}$
holds. $U^\prime$ is an open set of $U$. 

Further, we set 
$$\delta \ := \ 
\frac{1-\varepsilon^\prime}{1+\varepsilon^\prime}
\delta^\prime.$$
$\delta>0$ and 
$$\ \cup_{\lambda\in\Lambda^\prime} 
\, ( \, \{ \, y\in{\mathbb R}^m 
\, | \, \|y-f(0,\lambda)\|_{{\mathbb R}^m} 
< \delta \, \}
\times\{\lambda\} \, )
\ \subset \ V^\prime \ \subset \ 
{\mathbb R}^m\times\Lambda^\prime$$
hold. As $x\in {\mathbb R}^m$, $\|x\|_{{\mathbb R}^m}<\delta$, 
$\lambda\in \Lambda^\prime$ and $s\in[0,1]$ hold, 
because of $(sx,\lambda)\in G \subset E$, from $(12.1)$ and $(12.4)$, 
$\|f(x,\lambda)-f(0,\lambda)\|_{{\mathbb R}^m}
\leq(1+\varepsilon^\prime)\|x\|_{{\mathbb R}^m}
<(1+\varepsilon^\prime)\delta
=(1-\varepsilon^\prime)\delta^\prime,$
$(f(x,\lambda),\lambda) \in V^\prime$ and 
$(x,\lambda) \in U^\prime$ hold. Therefore, 
$$\{ \, x\in{\mathbb R}^m 
\, | \, \|x\|_{{\mathbb R}^m} < \delta \, \}
\times\Lambda^\prime
\ \subset \ U^\prime \ \subset \ 
{\mathbb R}^m\times\Lambda^\prime$$
holds. 

From $(12.1)$, 
$$\sup_{(x,\lambda)\in U^\prime}\|(D_xf)(x,\lambda)
-1_{{\mathbb R}^m}\|_{L({\mathbb R}^m;{\mathbb R}^m)}
\ \leq \ \varepsilon^\prime \ < \ \min\{\varepsilon,1\}$$
holds. From a usual inverse function theorem, 
for any $\lambda\in\Lambda^\prime$, the map 
$y \, \mapsto \, g(y,\lambda)$ is differentiable. Hence, 
for any $(y,\lambda)\in V^\prime$, 
$$(D_yg)(y,\lambda) \ = \ ((D_xf)(g(y,\lambda),\lambda))^{-1}$$
holds. The map $D_yg \, : \, V^\prime \, \rightarrow \, 
L({{\mathbb R}^m};{{\mathbb R}^m})$ is continuous. Also, 
$$\sup_{(y,\lambda)\in V^\prime}\|(D_yg)(y,\lambda)
-1_{{\mathbb R}^m}\|_{L({\mathbb R}^m;{\mathbb R}^m)}$$
$$= \ \sup_{(x,\lambda)\in U^\prime}\|((D_xf)(x,\lambda))^{-1}
-1_{{\mathbb R}^m}\|_{L({\mathbb R}^m;{\mathbb R}^m)}$$
$$= \ \sup_{(x,\lambda)\in U^\prime} \, \| \, 
\sum_{k=1}^\infty \, (1_{{\mathbb R}^m}-(D_xf)(x,\lambda))^k \, 
\|_{L({\mathbb R}^m;{\mathbb R}^m)}$$
$$\leq \ \sum_{k=1}^\infty \, {\varepsilon^\prime}^k
\ = \ \varepsilon$$
holds. 
\hfill 
$\blacksquare$

\newpage

{\bf References}

\vspace*{0.8em}

[1] B. W. Glickfeld, The Riemann sphere of a commutative Banach algebra, 
{\it Trans. Amer. Math. Soc.}, 134 (1968), 1-28. 

[2] K. R. Goodearl, Cancellation of low-rank vector bundles, 
{\it Pacific J. Math.}, 113 (1984), 289-302. 

[3] E. Hille and R. S. Phillips, {\it Functional analysis and semi-groups}, 
American Mathematical Society, 1957. 

[4] S. Kobayashi, Manifolds over function algebras and mapping spaces, 
{\it Tohoku Math. J.}, 41 (1989), 263-282. 

[5] E. R. Lorch, The theory of analytic functions in normed Abelian vector rings, 
{\it Trans. Amer. Math. Soc.}, 54 (1943), 414-425. 

[6] P. Manoharan, A nonlinear version of Swan's theorem, 
{\it Math. Z.}, 209 (1992), 467-479. 

[7] P. Manoharan, Generalized Swan's theorem and its application, 
{\it Proc. Amer. Math. Soc.}, 123 (1995), 3219-3223. 

[8] P. Manoharan, A characterization for spaces of sections, 
{\it Proc. Amer. Math. Soc.}, 126 (1998), 1205-1210. 

[9] A. S. Morye, Note on the Serre-Swan theorem, 
{\it Math. Nachr.}, 286 (2013), 272-278. 

[10] J. Mujica, {\it Complex analysis in Banach spaces}, North-Holland, 1986. 

[11] R. G. Swan, Vector bundles and projective modules, 
{\it Trans. Amer. Math. Soc.}, 105 (1962), 264-277. 

[12] L. N. Vaserstein, Vector bundles and projective modules, 
{\it Trans. Amer. Math. Soc.}, 294 (1986), 749-755.

\vfill

\noindent
Acknowledgment: 

\noindent 
This work was supported by JSPS KAKENHI Grant Number JP16K05245.

%\noindent Running title: Expanding fronts in an anisotropic diffusion equation  

%\vspace*{2.4em} 

\end{document}